\documentclass[12pt]{amsart}
\newcommand{\pile}{\pi^e_{ _\dw}}
\newcommand{\pilc}{\pi^c_{ _\uw}}

\newcommand{\res}{{\mtr{Res}}}

\newcommand{\epsm}{\eps_M}

\newcommand{\idm}{\id_M}

\newcommand{\zom}{{Z(M)}}
\newcommand{\csl}{{\cs_{\ml}}}

\newcommand{\mtcj}{\mtc J}

\newcommand\numeq[1]%
  {\stackrel{\scriptscriptstyle(\mkern-1.5mu#1\mkern-1.5mu)}{=}}

\newcounter{relctr} 
\everydisplay\expandafter{\the\everydisplay\setcounter{relctr}{0}} 

\newcommand{\cfl}{{\mtr{CF}}( L)}
\newcommand{\cel}{{\mtr{CE}}( L)}

\newcommand{\zresal}{\mtr{Res}^{ A}_{ L}}
\newcommand{\resal}{{\mtr{Res}}}

\newcommand{\zim}{\{0,\dots , m\}}
\newcommand{\zir}{\{0,\dots , r\}}

\usepackage{tikz-cd}
\usepackage{latexsym,amssymb,amsmath}
\usepackage{hyperref}
\usepackage{comment}
\usepackage{mathrsfs}
\usepackage{amsmath,amscd}
\usepackage[enableskew]{youngtab}
\usepackage[all]{xy}
\usepackage{empheq}

\usepackage{amsmath}

\usepackage{empheq}
\usepackage{xcolor}
\definecolor{lightgreen}{HTML}{90EE90}

\newcommand{\cocc}{\co(\cc)}

\newcommand{\mtj}{\mtc L}
\newcommand{\mtjcd}{\mtj_{\cd}}
\newcommand{\lag}{\langle}
\newcommand{\rag}{\rangle}

\newcommand{\evx}{{\ev_X}}

\usepackage[utf8]{inputenc}
\usepackage[T1]{fontenc}

 \newcommand{\vsk}{\vskip 0.15cm \noindent}

 \newcommand{\sumstor}{\sum_{s=0}^r}
 \newcommand{\fpcc}{\fp(\cc)}

\newcommand{\bxt}{\boxtimes}

\newcommand{\vect}{\mtr{Vec}}

\newcommand{\sent}{\mapsto}

\newcommand{\cec}{\mtr{  CE}(\cc)}

\DeclareMathOperator{\rev}{rev}

\newcommand\C{\mathcal{C}}
\DeclareMathOperator{\id}{id}

\DeclareMathOperator{\ev}{ev}
\DeclareMathOperator{\coev}{coev}

\excludecomment{verlong}
\includecomment{vershort}
\excludecomment{noncompile}

\usepackage{tikz}
\usetikzlibrary{matrix}
\usepackage{combelow}
\newcommand{\ccb}{{\mathcal B}}
 \newcommand{\bdelta}{\Delta}
 
\usepackage{amsmath,amsthm,amsfonts,amssymb}
\usepackage{amsmath, amsthm, amssymb, amscd, amsfonts}
\usepackage[all]{xy}
\usepackage{amssymb, amsthm, amsmath}
\usepackage{amsfonts}
\usepackage{color}
\DeclareMathOperator{\ad}{ad} 

\newcommand{\ra}{\rightarrow}

\newcommand{\ot}{\otimes}
\newcommand{\co}{\mathcal O}
\newcommand{\xra}{\xrightarrow}
\newcommand{\mtc}{\mathcal}
\newcommand{\cs}{\mtc S}

\newcommand{\lam}{\lambda}

\newcommand{\Lam}{\Lambda}

\newcommand{\al}{\alpha}
\newcommand{\eps}{\epsilon}

\newcommand{\ul}{\underline}

\newcommand{\lh}{\leftharpoonup}

\newtheorem{lem}[equation]{Lemma}

\theoremstyle{plain}

\newtheorem{prop}[equation]{Proposition}
\newtheorem{cor}[equation]{Corollary}
\newtheorem{rem}[equation]{Remark}

\newtheorem{thm}{Theorem}[section]
\numberwithin{equation}{section}

\newtheorem{defn}[thm]{Definition}

\newcommand{\dw}{\downarrow}
\newcommand{\uw}{\uparrow}

\newcommand{\ch}{\chi}
\newcommand{\mtr}{\mathrm}

\newcommand{\ncm}{\newcommand}
\ncm{\np}{\newpage}
\ncm{\ebl}{\end{thebibliography}}
\ncm{\bbl}{\begin{thebibliography}}
\ncm{\chd}{_{ _{\ch}}}
\ncm{\ald}{_{ _{\al}}}
\newcommand{\blam}{\Lam}
\ncm{\cP}{\mathcal{P}}
\ncm{\ei}{e_i}
\ncm{\eij}{e_{i,\;j}}
\ncm{\bt}{\begin{thm}}
\ncm{\bdef}{\begin{defn}}
\ncm{\edf}{\end{defn}}
\ncm{\et}{\end{thm}}
\ncm{\bc}{\begin{cor}}
\ncm{\bl}{\begin{lem}}
\ncm{\el}{\end{lem}}
\ncm{\bpf}{\begin{proof}}
\ncm{\epf}{\end{proof}}
\ncm{\ec}{\end{cor}}
\ncm{\ord}{\mtr{ord}}
\ncm{\er}{\end{rem}}
\ncm{\br}{\begin{rem}}
\ncm{\bn}{\begin}

\ncm{\bp}{\begin{prop}}
\ncm{\ep}{\end{prop}}
\ncm{\bd}{
\begin{document}}
\ncm{\ed}{\end{document}}
\ncm{\beq}{\begin{equation}}
\ncm{\beqn}{\begin{equation*}}
\ncm{\eeq}{\end{equation}}
\ncm{\eeqn}{\end{equation*}}
\ncm{\bea}{\begin{eqnarray}}
\ncm{\eea}{\end{eqnarray}}
\ncm{\beanon}{\begin{eqnarray*}}
\ncm{\eeanon}{\end{eqnarray*}}\ncm{\ek}{\eps|_K}\ncm{\diez}{\#}
\ncm{\bwt}{\bowtie}
\ncm{\cC}{\mtc{C}}\ncm{\cc}{\mtc{C}}
\ncm{\cX}{\mtc{X}}
\ncm{\wt}{\widetilde}
\ncm{\sg}{\sigma}
\ncm{\Rep}{\mathrm{Rep}}
\DeclareMathOperator{\Irr}{Irr}
\ncm{\X}{\mathcal{X}}
\ncm{\cA}{\mathcal{A}}
\ncm{\HKer}{\mtr{HKer}}
\ncm{\LKER}{\mtr{LKer}}
\ncm{\aad}{\mtr{ad}}
\newcommand{\mbf}{\mathbb F}
\ncm{\Dr}{\mtr{D}}
\ncm{\cD}{{\mathcal{D}}}\ncm{\cd}{{\mathcal{D}}}\ncm{\ce}{{\mathcal{E}}}
\ncm{\G}{\mathcal{G}}
\ncm{\Dc}{\mtc{D}}
\ncm{\E}{\mtc{E}}
\ncm{\fp}{\mtr{FPdim}}
\ncm{\Vc}{\mtr{Vec}}
\ncm{\cK}{\mtc{K}}
\ncm{\cM}{\mtc{M}}
\ncm{\cE}{\mtc{E}}
\ncm{\cS}{\mtc{S}}

\newcommand{{\ipr}}{i'}

\DeclareMathOperator{\End}{End}
\ncm{\cop}{\mtr{cop}}
\ncm{\op}{\mtr{op}}
\ncm{\chr}{character }\ncm{\ck}{\mtc{K}}
\ncm{\bw}{\bwt}
\ncm{\hker}{\mtr{HKer}}
\ncm{\bx}{\boxtimes}
\ncm{\blue}{\textcolor[rgb]{.00, .00, 1.00}}
\ncm{\red}{\textcolor[rgb]{1.00, .00, .00}}
\ncm{\green}{\textcolor[rgb]{.50, 0.20, .90}}
\ncm{\bne}{\begin{enumerate}}
\ncm{\ene}{\end{enumerate}}
\ncm{\lker}{\mtr{LKer}}
\ncm{\md}{\medbreak}
\ncm{\rep}{\Rep}\ncm{\ind}{\mtr{ind}}
\ncm{\mdn}{\md\noindent}
\ncm{\dd}{$}
\ncm{\up}{^}
\newcommand{\tcs}{\text}
\newcommand{\mbb}{\mathbb B}
\newcommand{\vs}{\mathbb V}
\newcommand{\sth}{suppose that\;}
\newcommand\rad{\operatorname{rad}}
\newcommand{\itm}{\item}
\newcommand{\dbd}{$$}
\newcommand{\mol}{\mtr{mod}}
 \newcommand{\ro}{\rho}
\newcommand{\irr}{\mathrm{Irr}}
\newcommand{\mbc}{\mathbb C}
\newcommand{\mbs}{\mathbb S}
\newcommand{\mbz}{\mathbb Z}
\newcommand{\ct}{\mtc T}
\newcommand{\sm}{\setminus}
\newcommand{\epl}{^{+}}
\newcommand{\sbsq}{\subseteq}
\newcommand{\sbs}{\subset}
\newcommand{\cco}{\mtr{co}}
\newcommand{\cz}{\mathcal{Z}}
\newcommand{\dual}{^{*}}
\newcommand{\Gm}{\Gamma}
\ncm{\cY}{\mtc{Y}}
\newcommand\ZZ{{\mathbb Z}} 
\newcommand{\bab}{\color{DarkOrchid}{}}
\newcommand{\eab}{\normalcolor{}}
\newcommand{\subs}{\subsection}
\newcommand{\cv}{\mtc{V}}
  \newcommand{\grn}{\green}
\newcommand{\dt}{\delta}

\newcommand{\ccf}{\mathrm{ {CF}(\cc)}}
\newcommand{\cce}{\mathrm{ {CE}(\cc)}}
\newcommand{\cecc}{\mathrm{ {CE}(\cc)}}
\newcommand{\cecd}{\mathrm{ {CE}(\cd)}}
\newcommand{\kk}{\Bbbk}
\newcommand{\otL}{\ot_{L}}
\newcommand{\otl}{\ot_{L}}
\newcommand{\unpsi}{1_{\psi}}
\newcommand{\epsi}{e_{\psi}}
\newcommand{\ephi}{e_{\phi}}
\newcommand{\ech}{e_{\ch}}
\newcommand{\nleftcid}{\text{left normal  coideal subalgebra}}
\newcommand{\dimL}{\dim_{\kk}L}
\newcommand{\cl}{\mtc L}
\newcommand{\mj}{\mtc J}
\newcommand{\tl}{\tilde L}
\newcommand{\tL}{\tilde L}
\newcommand{\tpsi}{\tilde(\psi)}
\newcommand{\tmx}{\tilde{\mtc X}}
\newcommand{\zlh}{\mathrm{ZL}}
\newcommand{\ba}{\mathrm A}
\newcommand{\bv}{\mathrm V}
\newcommand{\zhopf}{\mtc{Z}_{\mtr{Hopf}}}
\newcommand{\lstar}{L^{*}}
\newcommand{\ldstar}{L^{**}}
\newcommand{\mstar}{M^{*}}
\newcommand{\mdstar}{M^{**}}
\newcommand{\lkera}{\lker_{A}}
\newcommand{\mdprime}{M''}
\newcommand{\ldprime}{L''}
\newcommand{\cm}{\mtc M}
\newcommand{\ccm}{\mathcal M}
\newcommand{\cn}{\mathcal N}
\newcommand{\ccn}{\mathcal N}
\newcommand{\rx}{\mtr{Rex}}
\newcommand{\cca}{\ca}
\newcommand{\ih}{\underline{\mtr{Hom}}}
\newcommand{\cih}{\underline{\mtr{coHom}}}
\newcommand{\hm}{\mtr{ {Hom}}}
\newcommand{\cov}{\mtr{coev}}
\newcommand{\rora}{\rho^{\mtr{ra}}}
\newcommand{\rola}{\rho^{\mtr{la}}}
\newcommand{\cx}{\mtc X}
 \newcommand{\cZ}{\cz}
 \newcommand{\ca}{\cA}
 \newcommand{\stat}{\noindent}
 \newcommand{\bfa}{{\bf A}}
 \newcommand{\unu}{\mathbf{1}}
 \newcommand{\barzu}{{\bar {  Z}(\unu)}}
 
\newcommand{\idx}{\id_X}
\newcommand{\lprime}{L'}
\newcommand{\mprime}{M'}
\newcommand{\nat}{ \mtr{{  Nat}}}
\newcommand{\ft}{\mtc F_\lam}
\newcommand{\rhau}{\rightharpoonup}
\newcommand{\lhau}{\leftharpoonup}
\newcommand{\cf}{\mathrm{ {CF}}}

\newcommand{\cfc}{\mathrm{{CF}}(\cc)}
\newcommand{\csu}{{\mtr{C}}}
\newcommand{\cfcc}{\mathrm{ {CF}}(\cc)}
\newcommand{\cfcd}{\mathrm{CF}(\cd)}
\newcommand{\cfd}{\mathrm{CF}(\cd)}
\newcommand{\czcc}{{\cz(\cc)}}
\newcommand{\czcd}{{\cz(\cd)}}
\newcommand{\czt}{{\cz(\cz(\cc))}}
\newcommand{\enx}{\mtr{  End}}
\newcommand{\runu}{R(\unu)}

\newcommand{\bdfn}{\bn{defn}}
\newcommand{\edfn}{\end{defn}}
\newcommand{\deltax}{\delta_X}
\newcommand{\deltav}{\delta_V}
\newcommand{\repcca}{\rep_\cc(A)}
\newcommand{\xotay}{X \ot_A Y}
\newcommand{\xoty}{X \ot Y}
\newcommand{\votw}{V \ot W}
\newcommand{\votaw}{V \ot_A W}
\newcommand{\dimax}{\dim_AX}
\newcommand{\dimccx}{\dim_\cc(X)}
\newcommand{\dimcca}{\dim_\cc(A)}
\newcommand{\dimccv}{\dim_\cc(V)}
\newcommand{\dima}{\dim_A}
\newcommand{\biga}{A}
\newcommand{\comp}{\mathbb C}
\newcommand{\tehtaa}{\theta_A}
\newcommand{\tetaa}{\theta_A}
\newcommand{\ida}{\id_A}
\newcommand{\hma}{\hm_A}
\newcommand{\hmcc}{\hm_\cc}
\newcommand{\fv}{F(V)}
\newcommand{\fw}{F(W)}
\newcommand{\ota}{\ot_A}
\newcommand{\repza}{\rep_\cc^0(A)}
\newcommand{\epsa}{\eps_A}
\newcommand{\bndefn}{\bn{defn}}
\newcommand{\edefn}{\end{defn}}
\newcommand{\bdefn}{\bn{defn}}

\newcommand{\vld}{V^{*}}
\newcommand{\vldd}{V^{**}}
\newcommand{\xld}{X^{*}}
\newcommand{\xldd}{X^{**}}
\newcommand{\yld}{Y^{*}}
\newcommand{\yldd}{Y^{**}}
\newcommand{\aldu}{A^{*}}
\newcommand{\aldd}{A^{**}}

\newcommand{\ia}{\mtr{i}_A}
\newcommand{\aota}{A\ot A}

\newcommand{\idv}{\id_V}

\newcommand{\ld}{^*}
\newcommand{\repg}{\rep(G)}

\newcommand{\thetav}{\theta_V}

\newcommand{\tta}{\theta_A}

\newcommand{\muv}{\mu_V}
\newcommand{\muw}{\mu_W}

\newcommand{\dimcc}{\dim(\cc)}
\newcommand{\chii}{\chi_i}
\newcommand{\chistar}{\ch_{i^*}}
\newcommand{\chj}{\ch_j}
\newcommand{\chm}{\ch_m}
\newcommand{\chn}{\ch_n}
\newcommand{\dimvi}{\dim(V_i)}
\newcommand{\mtcd}{Q}
\newcommand{\mtca}{\mtc A}
\newcommand{\lamcd}{\lam_\cd}
\newcommand{\fpdimcd}{\fp(\cd)}
\newcommand{\laml}{\lam_L}
\newcommand{\apm}{A//M}
\newcommand{\apl}{A//L}
\newcommand{\repapm}{\rep(\apm)}
\newcommand{\repapl}{\rep(\apl)}
\newcommand{\dimvj}{\dim(V_j)}
\newcommand{\dvi}{\dim(V_i)}
\newcommand{\dvj}{\dim(V_j)}
\newcommand{\sumjtom}{\sum_{j=0}^m}
\newcommand{\sumitom}{\sum_{i=0}^m}
\newcommand{\sij}{s_{ij}}
\newcommand{\sji}{s_{ji}}
\newcommand{\dxj}{d_j}
\newcommand{\dxi}{d_i}
\newcommand{\dimka}{\dim_{\kk}(A)}
\newcommand{\dimk}{\dim_{\kk}}
\newcommand{\blaml}{\blam_L}
\newcommand{\sumjtor}{\sum_{j=0}^r}
\newcommand{\dimkl}{\dim_{\kk}(L)}
\newcommand{\mtcjl}{\mtc J_{ L}}
\newcommand{\vota}{ V\ot A}
\newcommand{\vi}{V_i}
\newcommand{\vj}{V_j}
\newcommand{\dimcd}{\dim(\cd)}

\newcommand{\alij}{\al_{ij}}
\newcommand{\alji}{\al_{ji}}
\newcommand{\rcc}{r_\cc}
\newcommand{\rcd}{r_\cd}
\newcommand{\clsx}{[X]}
\newcommand{\clsy}{[Y]}
\newcommand{\clsz}{[Z]}
\newcommand{\rcdp}{r_{\cd'}}
\newcommand{\sumjtorp}{\sum_{j=0}^{r'}}
\newcommand{\aljm}{\al_{jm}}
\newcommand{\aljn}{\al_{jn}}
\newcommand{\sjm}{s_{jm}}
\newcommand{\smj}{s_{mj}}
\newcommand{\snj}{s_{nj}}

\newcommand{\betaij}{\beta_{ij}}
\newcommand{\betaji}{\beta_{ji}}

 \newcommand{\ip}{{i'}}
\newcommand{\sumjtoprp}{\sum_{j=0}^{r'}}
\newcommand{\sumjtopr}{\sum_{j=0}^{r}}
 \newcommand{\teh}{\tilde{h}}
\newcommand{\cdp}{\cd'}
\newcommand{\xphii}{X_{\phi(i)}}
\newcommand{\inv}{^{-1}}

\newcommand{\fq}{f_Q}
\newcommand{\tr}{\mtr{tr}}
\newcommand{\rtwone}{R_{21}R}

\newcommand{\ccad}{{\cc_{\mtr{ad}}}}
\newcommand{\ccpt}{{\cc_{\mtr{pt}}}}
\newcommand{\qtr}{quasi-triangular\;}
\newcommand{\trq}{\tr_q}

\newcommand{\repal}{\mtr{Rep}(A//L)}
\newcommand{\lkeravi}{\lker_A(V_i)}
\newcommand{\lkeravj}{\lker_A(V_j)}
\newcommand{\cross}[1][1pt]{\ooalign{%
 \rule[1ex]{1ex}{#1}\cr
 \hss\rule{#1}{.7em}\hss\cr}}
\newcommand{\blml}{\blam_L} 
\newcommand{\phir}{\phi_R}
\newcommand{\kda}{{  \Phi(A)}}

\newcommand{\mtcil}{\mtc{I}_L}

\newcommand{\un}{\unu}
\newcommand{\tfl}{\mtc{T}}
\newcommand{\barzm}{\barz(M)}
\newcommand{\barzn}{\barz(N)}
\newcommand{\ccr}{\mtc R^{\cc}}
\newcommand{\ulc}{\ul{\cc}}

\newcommand{\pimx}{\pi_{M;\;X}}
\newcommand{\pinx}{\pi_{N;\;X}}
\newcommand{\acc}{{\mathrm A_\cc}}
\newcommand{\epsu}{{\bf \eps}_\unu}

\newcommand{\ob}{\mtr{Obj}}
\newcommand{\obc}{\mtr{Obj(\cc)}}
\newcommand{\ccop}{\cc^{\mtr{op}}}
\newcommand{\mtf}{\mtc F_\lam}
\newcommand{\mtfi}{\mtc F^{-1}_\lam}
\newcommand{\elcd}{\ell_\cd}
\newcommand{\mcid}{\mtc I_\cd}
\newcommand{\mcidp}{\mtc I_{\cd'}}
\newcommand{\wtildelcd}{\widetilde{\elcd}}
\newcommand{\wtildelcdp}{\widetilde{\ell_{\cd'}}}
\newcommand{\cpt}{\cc_{\mtr{pt}}}
\newcommand{\barzr}{\barz_\cd}
\newcommand{\barzv}{\barz(V)}
\newcommand{\acd}{\mathrm A_\cd}
\newcommand{\czrcd}{\cz_\cc(\cd)}
\newcommand{\sml}{\Small}
\newcommand{\bs}{\blue{\Small }}
\newcommand{\yd}{Yetter-Drinfeld}

\newcommand{\sumitor}{\sum_{i=0}^r}
\newcommand{\cdop}{\cd^{\mtr{op}}}
\newcommand{\ccrev}{\cc^{\mtr{rev}}}
\newcommand{\barz}{{\bar{\mathrm Z}}}
\newcommand{\etl}{\'etale\;}
\newcommand{\czca}{\cz(\ca)}

\usepackage{tikz-cd}
\usepackage{tikz-cd}
\usepackage{latexsym,amssymb,amsmath}
\usepackage{hyperref}
\usepackage{comment}
\usepackage{mathrsfs}
\usepackage{amsmath,amscd}
\usepackage[enableskew]{youngtab}
\usepackage[all]{xy}
\usepackage{empheq}
\newcommand{\diml}{\dim(L)}

\bd
\newcommand{\lah}{\;^A\lh}
\newcommand{\szccj}{|\cc^j|}
\newcommand{\bccm}{\mtr C_m}
\newcommand{\bccj}{\mtr C_j}
\newcommand{\bccl}{\mtr C_l}
\newcommand{\bccp}{\mtr C_p}
\newcommand{\zr}{{Z_r}}
\newcommand{\hpl}{{Z_r(\unu)}}
\newcommand{\rephpl}{\rep(H//L)}
\newcommand{\elm}{{\ell}}
\newcommand{\resl}{\res^A_L}
\newcommand{\mtrcj}{\mtr C_j}
\newcommand{\ccl}{{\cc_\ml^{\mtr{triv}}}}
\newcommand{\dtdots}{:}

\newcommand{\inccbs}{{i\in \ccb_s}}
\newcommand{\jinmtcjl}{{j\in \mtcjl}}
\newcommand{\fj}{F^j}
\newcommand{\iotael}{\iota^e}
\newcommand{\dimccj}{\dim(\cc^j)}
\newcommand{\ncj}{{\frac{\cc^j}{\dimccj}}}
\newcommand{\ma}{{\mtr A}}

\newcommand{\phiainv}{\phia^{-1}}
\newcommand{\elmstar}{{\eta_s}}
\newcommand{\dimccl}{\dim(\cc_L)}
\newcommand{\mtce}{\mtc E}
\newcommand{\ccmip}{\cc^{ ^{M(i')}}}

\newcommand{\zunu}{Z(\unu)}
\newcommand{\nzm}{{Z(M)}}
\newcommand{\ru}{R(\unu)}
\newcommand{\ml}{{\mtr L}}

\newcommand{\dimcsl}{\dim(\cs_\ml)}
\newcommand{\betal}{\beta(\ml)}
\newcommand{\alczml}{\al_{Z(M)}\big|_{ _\ml}}
\newcommand{\alcml}{{\al_M\big|_{ _\ml}}}
\newcommand{\mm}{{\mtr M}}
\newcommand{\mthn}{\mathbb N}

\newcommand{\csm}{\cs_\mm}
\newcommand{\mbp}{{\mtr P}}

\newcommand{\coevx}{\mtr{coev}_X}
\newcommand{\revx}{\mtr{ev}^r_X}
\newcommand{\rcoevx}{\mtr{coev}^r_X}
\newcommand{\xstar}{{X^*}}

\newcommand{\ccpj}{{\cc^{(j)}}}
\newcommand{\cfa}{\cfcc}
\newcommand{\cflo}{\cf(L_1)}
\newcommand{\cfltw}{\cf(L_2)}
\newcommand{\fjst}{F^{j}_{st}}
\newcommand{\ccpsj}{{\cc^{(j)}_s}}
\newcommand{\ccpjs}{\ccpsj}
\newcommand{\ccjps}{\ccpsj}
\newcommand{\ccjsp}{\ccpsj}
\newcommand{\tfjst}{\widetilde{\fjst}}
\newcommand{\lj}{\mtc L_j}
\newcommand{\mtmj}{\mtc M_j}

\newcommand{\tfjpsptp}{{\widetilde{F^{j'}_{s't'}}}}
\newcommand{\tfjstp}{{\widetilde{F^{j}_{st'}}}}
\newcommand{\enxczcca}{{\enx_\czcc(A)}}
\newcommand{\ccjpt}{\cc^{(j)}_t}
\newcommand{\ccpjt}{\cc^{(j)}_t}
\newcommand{\sumjtomst}{\sumjtom\;\sum_{s,t\in \mtmj}}
\newcommand{\fjpsptp}{{F^{j'}_{s't'}}}
\newcommand{\fjstp}{{F^{j}_{st'}}}
\newcommand{\alijst}{{\al(i)^j_{st}}}
\newcommand{\csujst}{\csu^j_{st}}
\newcommand{\st}{_{st}}
\newcommand{\sumrtom}{\sum_{r=0}^m}
\newcommand{\dimvr}{\dim(V_r)}
\newcommand{\fjss}{F^j_{ss}}
\newcommand{\csujss}{\mtr{C}^j_{ss}}
\newcommand{\csujts}{\mtr{C}^j_{ts}}
\newcommand{\jp}{{j'}}
\newcommand{\tp}{{t'}}
\newcommand{\spr}{{s'}}
\newcommand{\fjts}{F^j_{ts}}
\newcommand{\iotal}{\iota_{L}}

\newcommand{\hmcd}{{\hm_\cd}}
\newcommand{\csuj}{\csu^j}
\newcommand{\mujst}{\mu^j_{st}}
\newcommand{\mukuv}{\mu^k_{uv}}
\newcommand{\mujpsptp}{\mu^{j'}_{s't'}}
\newcommand{\cjst}{\csu^j_{st}}
\newcommand{\cjts}{\csu^j_{ts}}
\newcommand{\cjpsptp}{\csu^{j'}_{s't'}}
\newcommand{\cjptpsp}{\csu^{j'}_{t's'}}
\newcommand{\ncjst}{\frac{\cjst}{\dim(\cc_j)}}
\newcommand{\ncjts}{\frac{\cjts}{\dim(\cc_j)}}
\newcommand{\ncjpsptp}{\frac{\cjpsptp}{\dim(\cc_{j'})}}
\newcommand{\ncjptpsp}{\frac{\cjptpsp}{\dim(\cc_{j'})}}
\newcommand{\wcfcc}{{\widehat{\cfc}}}

\newcommand{\phia}{{\psi_{\ma,\unu}}}
\newcommand{\phil}{{\psi_{\ml,\unu}}}
\newcommand{\fpcd}{\fp(\cd)}
\newcommand{\wtd}{\widetilde}
\newcommand{\iotacl}{{\iota^c_{L}}}

 \newcommand{\fcjps}{{F(\ccjps)}}
\newcommand{\fcjpt}{{F(\ccjpt)}}
\newcommand{\mta}{\ma}
\newcommand{\otmtha}{{\ot_\ma}}
\newcommand{\mmino}{{m_{-1}}}
\newcommand{\nmino}{{n_{-1}}}
\newcommand{\bigop}{\big(}
\newcommand{\bigcp}{\big)}
 \newcommand{\dmino}{{_{-1}}}
  \newcommand{\dzero}{{_{0}}}
\newcommand{\cresal}{\mtr{CFRes}^A_L}
\newcommand{\pil}{{\pi}}

\newcommand{\epsul}{{\ul{\eps}_L}}
\newcommand{\acml}{\al_M\big|_L}
\newcommand{\racml}{{\al_M\big|_{ _L}}}
\newcommand{\betam}{{\beta( M)}}
\newcommand{\twd}{:}
\newcommand{\actcm}{{\al_M}}
\newcommand{\rsimengv}{{\simeq_{\mtr{REGNV}}}}
\newcommand{\deltaunu}{\delta_{\unu}}
\newcommand{\dedots}{:}
\newcommand{\lamz}{{{\lam_0}}}

\newcommand{\mtcls}{{{\mathcal L}_s}}

\newcommand{\dimca}{\dim(A)}

\newcommand{\dimrs}{\dim(r_s)}
\newcommand{\ccbs}{{\ccb_s}}
\newcommand{\brs}{{R_s}}
\newcommand{\dimbrs}{{\dim(R_s)}}
\newcommand{\fst}{F_{st}}
\newcommand{\enxczcc}{{\enx_{\czcc}}}
\newcommand{\ellcd}{{\ell_{\cd}}}
\newcommand{\etas}{\eta_s}

\newcommand{\tinl}{\blue{{\ul{t} \in \mtc L}}}
\newcommand{\almcd}{\lamcd}
\newcommand{\idual}{{i^*}}
\newcommand{\ipdual}{{{\ip}^*}}
\newcommand{\deltast}{\delta_{s,t}}
\newcommand{\dip}{{d_\ip}}
\newcommand{\mtfinv}{{\mtf^{-1}}}
\newcommand{\wellcd}{{\widetilde{\ell_\cd}}}
\newcommand{\pie}{\pile}
\newcommand{\sinl}{{s\in \lj}}
\newcommand{\bls}{{s}}
\title[Fusion categories]{Subalgebras of \'etale algebras in braided fusion categories}

\author{Sebastian Burciu}
\address{Inst.\ of Math.\ ``Simion Stoillow" of the Romanian Academy
P.O. Box 1-764, RO-014700, Bucharest, Romania}
\email{sebastian.burciu@imar.ro}
\thanks{This work was supported by a grant of the Ministry of Research, Innovation and Digitization, CNCS/CCCDI - UEFISCDI, project number PN-III-P4-ID-PCE-2020-0878, within PNCDI III}
\date{\today}
\maketitle
\begin{abstract}

 In \cite[Rem. 3.4]{DNO} the authors asked the question if any \'etale subalgebra of an \'etale algebra in a braided fusion category is also \'etale. We give a positive answer to this question if the braided fusion category $\cc$ is pseudo-unitary and non-degenerate.   
In the case of a pseudo-unitary fusion category we also give a new description of the lattice correspondence from \cite[Theorem 4.10]{DMNO}. This new description enables us to describe the two binary operations on the lattice of  fusion subcategories.
\end{abstract}

\section{Introduction}
In \cite[Rem. 3.4]{DNO} the authors asked the question if any \'etale subalgebra of an \'etale algebra in a braided fusion category $\cc$ is also \'etale. It seems that an answer  to this question does not appear in
the literature at this moment. We give a positive answer to this question if the braided fusion category $\cc$ is pseudo-unitary and non-degenerate.

For any fusion category $\cc$ it was shown in \cite[Lemma 3.5]{DMNO} that $\ma=R(\unu)$ is a connected \'etale algebra in $\czcc$ where $R$ is a right adjoint of the forgetful functor $F:\czcc\ra \cc$. Moreover, in \cite[Theorem 4.10]{DMNO} it was shown that the lattice of fusion subcategories of $\cc$ is in a reverse one to one correspondence with the lattice of \'etale subalgebras of $\ma$. $\ma$ is sometimes called {\it the adjoint algebra} of $\cc$ in \cite{scalg}.

Recall that given a connected \'etale subalgebra $\ml$ of $\ma$  in \cite{DMNO} it was associated a fusion subcategory $\beta(\ml)$ of $\cc$. Under the tensor equivalence $R:\cc\ra \czcc_\ma$ the fusion subcategory $\betal$ consists of those $\ma$-modules in $\czcc$ that are dyslectic as right $\ml$-modules.

For a pivotal tensor categories, Shimizu has developed in \cite{scalg} a character theory for fusion categories, similar to the classical theory of representations of finite groups and semisimple Hopf algebras. Moreover he showed that the adjoint algebra $\ma$ associated to a fusion category has a natural action on any object of $\cc$.

The main goal of this paper is to show that in the case of a pseudo-unitary fusion category $\cc$ all unitary subalgebras of $\ma$ in $\czcc$ are connected \'etale. We also obtain a new characterization of the fusion subcategories from \cite{DMNO} associated to \'etale subalgebras of $\ma$.

In Section \ref{charthry}, in the case $\cc$ is a pseudo-unitary fusion category, we associate to any unitary subalgebra $\ml$ of $\ma$ in $\czcc$ a fusion subcategory $\cs_\ml$ of $\cc$, consisting of those simple objects of $\cc$ whose characters (as defined in \cite{scalg}) have trivial restriction to $\ml$. We also denote by ${\ccl}$ the full abelian subcategory of $\cc$ consisting of those objects that receive a trivial action from the subalgebra $\ml$, see Definition \ref{ccl}.

It is shown  that for  any pivotal fusion category $\cc$ one has a chain of inclusions of abelian full subcategories:
\beq\label{chain}
\beta(\ml)\subseteq {\ccl}\subseteq \cs_\ml.
\eeq

On the other hand, we show in Theorem \ref{fpdsl} that if $\cc$ is a pseudo-unitary fusion category then $\fp(\cs_\ml)=\frac{\fp(\cc)}{\fp(\ml)}$.
Since by \cite[Theorem 4.10]{DMNO} we also have $\fp(\betal)=\frac{\fp(\cc)}{\fp(\ml)}$ this shows that  the above chain inclusions of fusion subcategories are in fact equalities.  This chain of equalities implies that ${\ccl}$ is also a fusion subcategory and that any unitary connected subalgebra of $\ma$ is also an \'etale subalgebra of $\ma$. These facts  can be stated as follows:
\bt\label{main1}
Any fusion subcategory of  a pseudo-unitary fusion category $\cc$ is of the form ${\ccl}$ for some unitary subalgebra $\ml$ of the algebra $\ma:=R(\unu)$. Moreover, every such unitary subalgebra $\ml$ of $\ma$ is \'etale.
\et

From Theorem \ref{main1}, using some results concerning \'etale algebras developed in \cite{DMNO}  we deduce that in a non-degenerate braided pseudo-unitary fusion category, any subalgebra of an \'etale algebra is also \'etale. More precisely we prove the following:
\bt\label{main2}
Let $\cc$ be a pseudo-unitary non-degenerate braided fusion category and $A$ an \'etale algebra of $\ccb$. Then any unitary subalgebra $L$ of $A$ is also \'etale.
\et

For two subalgebras of the adjoint algebra $\ma$ we denote by $\ml\mm$ the image of $\ml\ot \mm$ under the multiplication $m:\ma\ot \ma\ra \ma$. Since $\ma$ is commutative it follows that $\ml\mm =\mm\ml$ and moreover $\mm \ml$ is a subalgebra of $\ma$. Concerning the two binary operations on the lattice of fusion subcategories we prove the following
\bt\label{main3}  
Let $\cc$ be a pseudo-unitary fusion category and $\ml,\mm$ two unitary subalgebras of the adjoint algebra $\ma$ of $\cc$. With the above notations one has
$$\csl\cap\csm=\cs_{\ml\mm}, \;\csl\vee \csm=\cs_{\ml\cap \mm}.$$
\et
Shortly, the organization of the paper is as follows. In Section \ref{prelim} we recall the basics on fusion categories and also the  character theory developed by Shimizu in \cite{scalg} for pivotal fusion categories. In Section \ref{charthry} we construct the fusion subcategory ${\cs_\ml}$ mentioned above in the pivotal case. We also define for any fusion category the full abelian subcategory $\ccl$ and prove the inclusion ${\cs_\ml}\supseteq \ccl$  from Equation \eqref{chain}. In Section \ref{coindm1} we prove the inclusion $\ccl\supseteq \betal$ from Equation \eqref{chain}. Also in this section we give the proofs of the above three theorems. 
In Subsection \ref{ha} we also consider as an example the case of fusion categories coming from representation categories semisimple Hopf algebras.


We work over an algebraically closed field $\kk$ of characteristic zero. All the Hopf algebras and fusion categories notations follow  \cite{EGNO15}.
\section{Preliminaries}\label{prelim}
In this section we review the basic properties of fusion categories that are needed through the paper. 
For the definition and  standard theory of monoidal categories, we refer the reader to \cite{ML98} and \cite{Kas}. Given a monoidal category we define by $\cocc$ the class of all its objects.
Recall that a {\it left dual object} of $X \in \cocc$ is a triple $(X^*, \evx, \coevx)$ consisting of an object $X^* \in \cocc$
and morphisms \dd \evx : X^* \otimes X \ra \unu$ and \dd \coevx : \unu \ra X \otimes X^*$ such that the following equalities are satisfied:
\beq\label{zz1r}
(\evx \otimes \id_{\xstar}) \circ (\id_{\xstar} \ot d) = \id_{\xstar},\;\; (\idx \ot \evx) \circ (\coevx \ot \idx) = \idx.
\eeq 
Similarly,  one can define a  {\it right dual} of $X$ (which in fact is a left dual of $X$ in $\cc^{rev}$). A monoidal category $\cc$ is said to be {\it rigid} if every object of $\cc$ has both a left and a right dual object. 

Recall that a {\it finite tensor category } \cite{EO} over a field $\kk$ is a rigid monoidal category $\cc$ such that $\cc$ is a finite abelian category,  the tensor product $\ot:\cc \times \cc \ra\cc$ is $\kk$-linear in each variable, and $\enx_\cc(\unu)\simeq\kk$ as algebras. A {\it fusion category} \cite{ENO} is a semisimple finite tensor category.

For a monoidal category $\cc$, the {\it left monoidal center} (or the Drinfeld center) of $\cc$ is a category $\czcc$ defined as follows: An object of $\czcc$ is a pair $(V, \sg_V)$ consisting of an object $V\in \cc$ and a natural isomorphism
$$\sg_{V, X}: V \ot  X \ra X \ot V $$
for all $X \in \cocc$, satisfying a part of the hexagon axiom. A morphism  $f:(V, \sg_V)\ra (W, \sg_W)$ in $\czcc$
 is a morphism in $\cc$  such that $(\id_X\ot f) \circ\sg_{V, X}=\sg_{W, X}\circ(f\ot \id_X)$ for all $X\in \cc$. The composition of morphisms is defined in an obvious way. The category $\czcc$ is in fact  a braided monoidal category, see, e.g., \cite[Chapt. XIII.3]{Kas} for details.

Let $\cc$ be any finite tensor category and $A$ be an algebra in $\cc$. Recall that the algebra $A$ is called {\it connected} if $\hm_\cc(\unu, A)=\kk$. A subalgebra $B\xra{\iota} A$ is called {\it unitary} if
$u_B\circ \iota=u_A$ where $u_B$ and $u_A$ are the units of $B$ and $A$.
Recall \cite{DMNO} that an algebra $A$ of a braided fusion category $\ccb$  is called {\it \'etale} if it is {\it separable} and connected. An algebra $A$ is called separable if its multiplication
$A\ot A\xra{m} A$ has a section as $A$-bimodules in $\cc$.
This is equivalent by \cite[Theorem 3.2]{DMNO} to the fact that the category of right (or left) $A$-modules $\ccb_A$ is a semisimple category.

Let $A$ be an algebra in a braided fusion category $\ccb$  and $M$ be a right $A$-module in $\ccb$. Then $M$ is called {\it dyslectic (local)} if $c_{A, M}\circ c_{M, A}\circ \ro_M=\ro_M$ where $\ro_M:M\ot A\ra M$ is the module structure of $M$.

Dyslectic modules form a full subcategory of $\ccb_A$ which is usually denoted by $\ccb^0_A$. This  subcategory is closed under $\ot_A$ and the braiding in $\ccb$ induces a natural braided structure in $\ccb^0_A$, see \cite[Section 2]{par}. Thus, $\ccb^0_A$
is a braided fusion category.
\subsection{The central Hopf comonad of a finite tensor category}\label{hcm} Let $\cc$ be a fusion category. The forgetful functor $F:\czcc\ra \cc$ admits a right adjoint functor $R:\cc \ra \czcc$  and $Z :=FR:\cc \ra \cc$ is a Hopf comonad called {\it the central Hopf comonad associated} to $\cc$.  Moreover, one has that 
\beq
Z(V)\simeq \int_{X\in \cc}X\ot V\ot X^*.
\eeq
see \cite[Section 2.6]{scalg}. We denote by $\pi_{V;X}:Z (V)\ra X\ot V\ot X^*$ the universal dinatural transformation associated to the end $Z (V)$.
The Hopf comonad structure of $Z $ can be described in terms of the dinatural transformation $\pi$. The comultiplication $\delta : Z  \ra {Z }^2$ is the unique natural transformation such that
\beq\label{cobz}
(\id_X \otimes \pi_{V ;Y}\otimes \id_{X^*} ) \circ \pi_{Z (V );X} \circ \delta_V = \pi_{V ;X\otimes Y}
\eeq
The counit of ${\bf {\eps}}:Z\ra \id_\cc $ is given by ${\bf \eps}_V:= \pi_{V ;1}.$
There is a lax monoidal structure defined on $Z$ by $$Z_2(M,N): Z(M)\ot Z(N)\ra Z(M\ot N)$$
as the unique map making the following diagram commutative:
\vsk
{\Small
\begin{tikzcd}[column sep=large,row sep=large]
 Z(M)\ot Z(N)
 \arrow[rrr,"Z_2{(M,N)}"] \arrow[d,"{\pi}_{M;\;X}\ot {\pi}_{N;\;X}"] 
 & 
 &
 & 
 Z(M\ot N) \arrow[d,"\pi_{M\ot N,\;X}"]
 \\ 
 (X\ot M \ot X^* ) \ot (X\ot N\ot X^*)
 \arrow[rrr,"\id_{X\ot M} \ot \ev_X\ot \id_{N\ot X^*}"]  & 
 &
 & 
 X\ot M\ot N \ot  X^*  
\end{tikzcd}
}

which can be written as
\beq\label{ztw}
\pi_{M\ot N,\;X}\circ Z_2(M,N)=(\id_{X\ot M} \ot \ev_X\ot \id_{N\ot X^*})\circ ({\pi}_{M;\;X}\ot {\pi}_{N;\;X})
\eeq

\subsection{On the adjoint algebra of a finite tensor category}
It is known that $\ma:=R(\unu)$ has the structure of central commutative algebra in $\cz(\cc)$. The half braiding of $\ma$, denoted by $c_{A, X}$, is defined by
$$
A\ot X\xra{\delta_\unu\ot \idx} Z(A)\ot X \xra{\pi_{Z(\unu),X}\ot \idx}X\ot A\ot X^*\ot X\xra{\idx\ot \id_A\ot \ev_X} X\ot A.
$$

The multiplication $m:\ma\ot \ma \ra \ma$ and the unit $u_\ma:\unu\ra \ma$ of the
adjoint algebra $\ma$ are uniquely determined by  by the universal property of the end $A=Z(\unu)$ as:
\beq\label{unit} \pi_{\unu;X} \circ u_\ma = \cov_X,
\eeq
\beq \label{ma}\pi_{\unu;X} \circ m = (\id_X \otimes \ev_X \otimes \id_{X^*} ) \circ (\pi_{\unu;X} \otimes \pi_{\unu;X}).
\eeq
Moreover $\epsu:\ma\ra \unu$ is a morphism of algebras, see \cite{scalg}. 

Any object $X\in \cc$ is canonically an $A$-module in $\cc$, via the morphism:
 \beq\label{rox}
\al_X:A\ot X\xra{\pi_{\unu; X}\ot \id_X}X\ot X^*\ot X\xra{\id_X\ot ev_X}X.
 \eeq
 We denote by $c_{R(M), X}:R(M)\ot X\ra X\ot R(M)$ the half braiding of $R(M)\in \czcc$. Recall that with the above braiding $\ma=R(\unu)$ becomes a commutative algebra in the center $\czcc$.
 \noindent
Moreover, by \cite[Equation (3.12)]{scalg} one has that
 \beq\label{rox2}
 \al_X=(A\ot X\xra{c_{\ma,X}}X\ot A\xra{\id_X\ot \eps_\unu}X).
 \eeq
The vector space $\cecc:= \hm_{\C}(\unu, A) $ is called {\it the set of central elements.} For $a, b \in \cecc$ we set $ab :=m \circ (a \ot b)$. Then the set $\cecc$ is a monoid with respect to this operation. Note that for the  algebra unit $u_\ma:\unu \ra \ma$ of the algebra $\ma$ to the unit of the monoid $\cecc$ is $F(u_\ma)$.
\subsection{Internal characters for pivotal strctures}
Recall that a pivotal structure $j$ on a tensor category $\cc$ is a tensor isomorphism $j:\id_\cc\ra (-)^{**}$.  Using the pivotal structure one can construct a {\it right evaluation} as follows:
{\Small
$$
\widetilde{ev}_X:X\ot X^*\xra{j\ot id}X^{**}\ot X^*\xra{ev_{X^*}} \unu.
$$
} 
Then {\it the right partial pivotal trace} of $f:A\ot X\ra B\ot X$ is defined as follows:
{\Small
\beq\label{ptr}
\tr_{A, B}^X:A=A\ot \unu\xra{\id_A \ot coev_X}A\ot X \ot X^*\xra{f \ot id} B\ot X\ot X^*\xra{\id_B  \ot \widetilde{ev}_X} B.
\eeq
}
The usual {\it right pivotal trace} of an endomorphism $f :X\ra X$ is obtained as a particular case for $A=B=\unu$. In particular,  the {\it right pivotal dimension $\dim^r(X)$ of $X$} is defined as the right trace of the identity of $X$. A pivotal structure a on a tensor category $\cc$ is called {\it spherical} if
$\dim(V)=\dim(V^*)$ for any object $V \in \co(\cc)$. 
Given an object $X \in \co(\cc)$ the internal character $\mtr{ch}(X)$ is defined as the  morphism 
{\Small
\beqn
\mtr{ch}(X):=\tr^{X}_{A, \unu}(\al_X):A\ra \unu.
\eeqn
}
Therefore using Equation \eqref{ptr} one can write that
{\Small
\beq\label{chx}
\mtr{ch}(X)=A\ot \unu\xra{\id \ot coev_X}A\ot X \ot X^*\xra{\al_X\ot id} X\ot X^*\xra{\widetilde{ev}_X} \unu.
\eeq
}
Then the space $\cfcc:=\hm_\cc(A, \unu)$ is called the {\it space of class functions} of $\cc$. For two class functions $f, g\in \cfcc$ one can define a multiplication by
 $$f\star g:=f \circ Z(g) \circ \delta_{\unu}.$$ 
 Here $\delta: Z \ra  Z^2$ is the comultiplication structure of $Z$ defined in the Equation \eqref{cobz}. By \cite[Theorem 3.10]{scalg} one has that $\mtr{ch}(X\ot Y)=\mtr{ch}(X)\mtr{ch}(Y)$ for any two objects $X$ and $Y$ of $\cc$. 
For  a finite tensor category  $\cc$ the space of class functions $\cfcc$ is a finite-dimensional algebra.
\subsection{The cointegral and integral of a fusion category}
Let $\cc$ be a fusion category and $A=Z(\unu)$ be its adjoint algebra as defined above.

 An {\it integral} in $\cc$ is a morphism $\Lambda: \unu \ra A$   in $\cc$   such that
 {\Small
 \beqn
 m \circ (\id_{A}\otimes \Lambda)=\epsu\ot \Lam.
 \eeqn
 }
\noindent
A {\it cointegral} in $\cc$  is a morphism $\lambda : A\ra \unu $  such that
{\Small
\beqn
Z(\lam)\circ \delta_{\unu}=u \otimes \lam
\eeqn
}
where $u:\unu \ra A$ is the unit of the algebra $A$. It is well known that the integral and cointegral of a finite unimodular tensor category are unique up to a scalar, see \cite{scalg}.
\subsection{Fourier transform for finite tensor categories}\label{ft}
Let  $\lambda\in \cfcc$ be  a non-zero integral of $\cc$. The {\it Fourier transform for finite tensor categories} of $\cc$  associated to $\lambda$ is the linear map
\beq
\mtc F_{\lambda}:\cecc\ra \cfcc\;\;\text{given by}\;\;a \mapsto \lambda \lh \mtc S(a)
\eeq

\noindent
where $\mtc S: \cecc\ra \cecc$ is the antipodal operator on $\cecc$, see \cite[Definition 3.6]{scalg}. The Fourier transform is a bijective $\kk$-linear map whose inverse is given in  \cite[Equation (5.17)]{scalg}. Here, the right action denoted by $\lh$ of $\cecc$ on $\cfcc$ given by $f \lh b=f \circ m \circ (b\ot \id_{A})$ for all $f \in \cfcc$ and $b \in \cecc$. There is also a non-degenerate evaluation pairing $\langle,\;\rangle_{\mtr{\ma} }$ given by
\beq\label{eva}
\langle\;,\;\rangle_{\mtr{\ma}} :\cfcc \ot \cecc  \ra \kk,\;\langle \ch, z\rangle\sent \ch\circ z.
\eeq

\subsubsection*{Fourier transform for fusion categories}
For the rest of this section, suppose that $\cc$ is a pivotal fusion category over an algebraically closed field $\kk$. Furthermore, let $\irr(\cc)\twd=\{V_0, \dots, V_m\}$ be a complete set of representatives of isomorphism classes of simple objects with $V_0=\unu$, the unit object. For $i \in \{0, \dots, m\}$, we define $i^*\in \{0, \dots, m\}$ by $V_i^*\simeq V_{i^*}$. Then $i \sent i^*$ is an involution on $\{0, \dots, m\}$.
As an object of $\cc$, the adjoint algebra decomposes as
\beq
A\simeq \bigoplus_{i=0}^m V_i\ot V_{i}^*.
\eeq
Shimizu has defined in  \cite{scalg}  the elements
$$
E_i:\unu \xra{\mtr{coev}_{V_i}} V_i\ot V_i^* \hookrightarrow A,\;\; \ch_i: A \xra{\pi_i} V_i\ot V_i^* \xra{\widetilde{\mtr{ev}_{V_i}}} \unu.
$$
It is easy to see that $\{E_i\}_{i=0, \dots , m}$ and $\{\ch_i\}_{i=0, \dots , m}$ are bases for $\cecc$ and $\cfcc$ respectively, such that
\beq\label{fbases}
\lag\ch_i, E_j\rag=d_i\delta_{i, j}.
\eeq
where $d_i\twd=\dim(V_i)$. Moreover, $E_iE_j=\delta_{i,j}$ and $\mtc S(E_i)=E_{i^*}$, where $\tilde S:\cecc \ra \cecc$ is the antipodal mentioned above. The elements $\ch_i$ are called the irreducible characters of the simple objects $V_i$ and $E_i\in\cecc$ their corresponding primitive idempotents.
Note that $E_0$, the idempotent associated to the unit object $\unu$ of $\cc$, is the idempotent integral $\blam \in \cecc$, see \cite[Lemma 6.1]{scalg}. By \cite[Equation (6.8)]{scalg} one has that the idempotent cointegral of $\cc$ can be written as\twd
\beq\label{intregform}\lam_{\cc}=\frac{1}{\dimcc}(\sum_{[V_i]\in \irr(\cc)}d_{i^*}\ch_{i}).
\eeq 
Since $\blam=E_0$ it follows by Equation \eqref{intregform} that
$\langle\lam_\cc, \blam\rangle=\frac{1}{\dim(\cc)}$.
It also follows that
\beq\label{610}
\mtf^{-1}(\ch_i)=\frac{\dimcc}{d_i}E_{i^*}
\eeq
This equation was proven in \cite[Equation (6.10)]{scalg} under the hypothesis that the Grothendieck ring of $\cc$ is commutative but the proof works word by word also in the case of an arbitrary pivotal fusion category, with the  not necessarily commutative Grothendieck ring.

For a spherical fusion category $\cc$  one has $\dim(V_{i^*})=\dim(V_i)$ and it follows from  \cite[Equation (4.7)]{ccc-march} that$\twd$
\beq\label{comp1}
\langle \ch, \; \mtf^{-1}(\mu)\rangle_\ma =\dimcc \tau(\ch\mu).
\eeq
for any $\ch, \mu\in \cfcc$.
 \subsection{On the inclusion of class functions of a fusion subcategory}\label{cfinc}
Let $\cd$ be a fusion subcategory of a given fusion category $\cc$.  As in \cite[Sect. 4.3]{scalg},  for any object $V$ of $\cc$ the end 
$$
\bar Z(V):\cd^{\op}\times \cd\ra \cc,\;\; \bar Z(V):=\int_{X\in \cd}  X\ot V\ot X^*
$$ exists 
and we denote by $\bar{\pi}_{\unu\; X}:\barzu \ra X\ot X^*$ the universal dinatural maps defining this end. From the universal property of $\barzu$ there is a unique canonical map  $Z(V)\xra{q_V}\barz(V)$ in $\cc$, such that $\bar{\pi}_{V;\; X}\circ q={\pi}_{V;\; X}$ for any  object $X$ of $\cd$. In \cite[Appendix]{ccc-march} it is shown that $q:Z(\unu)\ra \barzu$ induces a map $q_\unu^*:\cfcd \ra \cfcc, \;\ch\sent \ch \circ q_\unu$
that  is a monomorphism of $\kk$-algebras.
\subsection{The adjunction isomorphisms for $F:\czcc\ra\cc$}
Let $\cc$ be any fusion category and $R:\czcc\ra \cc$ be a right adjoint of the forgetful functor as in the previous section. It is known that $Z:=FR$ is a Hopf comonad in $\cc$. One has that the category of $Z$-comodules in $\cc$ is equivalent to the center $\czcc$, see \cite{ds07} and \cite{bv12}. Moreover, in this case, the comodule forgetful functor $\cc^Z\ra \cc$ identifies to $F:\czcc\ra \cc$. This shows that the adjunction isomorphisms are given by 
\beq\label{can}
\psi_{V, \unu}:\hmcc(F(V), \unu)\ra \hm_{\czcc}(V, \ma) , \; f \sent  R(f)\circ \eta_V,
\eeq
and its inverse is given by 
\beq\label{caninv}
\psi_{V, \unu}^{-1}:\hm_{\czcc}(V, \ma) \ra\hmcc(F(V), \unu), \; g \sent \epsu\circ F(g).
\eeq
for any object $V\in \czcc$. Here $\eta_V:V\ra R(F(V))$ is the unit of the adjunction which coincides to the $Z$-comodule structure of $V\in \czcc$.
\subsection{Frobenius-Perron dimensions and pseudo-unitarity}
For a simple object $X$ of a fusion category $\cc$ the Frobenius-Perron eigenvalue of the left
multiplication by $[X]$ on the Grothendieck ring $K_0(\cc)$ is denoted by $\fp(X)$ and it is called the {\it Frobenius-Perron
dimension} of $X$. Recall that $\fp(X)$  is a positive real algebraic number and the Frobenius-Perron dimension extends linearly to an algebra morphism $\fp: K_0(\cc)\ot_{\mathbb Z}\comp\ra \comp$. The Frobenius-Perron dimension $\fp(\cc)$ of $\cc$ is defined as
$$
\fp(\cc):= \sum_{X\in \irr(\cc)}\fp(X)^2.
$$

A fusion category $\cc$ is called {\it pseudo-unitary} if $\fpcc =\dimcc$. In this case, $\cc$ admits a unique (canonical) spherical structure with respect to
which the categorical dimensions of simple objects are all positive, see \cite[Proposition 8.23]{ENO}. With respect to this spherical structure, the categorical dimension of any object coincides with its Frobenius-Perron dimension, i.e. $\fp(X)=\dim(X)$ for any object $X\in \co(\cc)$.  Moreover, every full fusion
subcategory of $\cc$ is pseudo-unitary.
\section{On the subcategory ${\cs_\ml}$ via character theory}\label{charthry}
Let as above $\cc$ be a pivotal fusion category and $F:\czcc\ra\cc$ be the forgetful functor with right adjoint $R:\cc\ra \czcc$. Let also $\ma:=R(\unu)$  be its adjoint algebra and $\iota:\ml\hookrightarrow \ma$ be a unitary subalgebra of $\ma$. We also denote by  $L$ and $A$ the images of $L$ and $A$ under the forgetful functor $F$. Since $\ma$ is connected it follows that $\ml$ is also connected.  Note also that $\ml$ is a commutative algebra in $\czcc$ as a subalgebra of $\ma$.  Since $\czcc$ is a fusion category there is also a projection  $\ma \xra{\pi} \ml$ in $\czcc$ of $\ma$ onto the object  $\ml$ such that $\pi \circ \iota=\id_\ml$. 
Since $\cc$ is a fusion category, the left adjoint of $F$ is also isomorphic to $R$ (see \cite{sh-imrn}), and therefore $\hm_\czcc(\unu, \ma)=\hm_\cc(\unu, \unu)\cong \kk$ . Thus the unit $u_\ml$ of the subalgebra $\ml$ can be written as $\pi\circ u_\ma$, where $\ma \xra{\pi} \ml$ is the projection in $\czcc$ of $\ma$ into $\ml$. 

Define also a {\it character space} of $\ml$ by $\cfl:=\hm_\cc(L, \unu)$. Then there are well-defined  restriction maps:
\beq
\resal:=\hm_\cc(F(\iotal)\text{,}\;-):\cfcc\ra\cfl,
\eeq
\beq
\zresal:=\hm_{\czcc}(\iotal\text{,}\;-):\hm_\czcc(\ma, \ma)\ra \hm_\czcc(\ml,\ma). 
\eeq
Since $Z=FR$ is a comonad it follows that the adjunction isomorphisms $F\dashv R$ are given in this case by 
\beq\label{adjisoml}
\psi_{\ml, \unu}:\hm_\cc(L,\unu) \xra{} \hm_{\czcc}(\ml, \ma),\;\; \mu\mapsto R (\mu)\circ \eta_\ml.
\eeq
where $\eta_\ml$ is the unit of the adjunction $L\xra{\eta_\ml} R(F(\ml))$. 

Moreover, its inverse is given by 
\beq\label{can}
\psi_{\ml, \unu}^{-1}:\hm_{\czcc}(\ml, \ma)\ra\hm_\cc(L,\unu), \; g \sent \epsu\circ g.
\eeq
Naturality in the first variable of the adjunction isomorphisms $\psi_{(-),\; \unu}$ of  $F\dashv R$ implies that the following diagram is commutative:
\begin{equation}\label{daj}
{
\begin{tikzcd}[column sep=5cm, row sep=1cm]
\cfl=\hmcc(L, \unu)
 \arrow[r,"\psi_{\ml, \unu}"]   
 &  
\hm_\czcc(\ml,\ma)  
  \\
\cfcc=\hm_\cc(A,\unu) \arrow[u, "\resal"] \arrow[r,"\psi_{\ma, \unu}"] 
&  
\hm_\czcc(\ma, \ma) \arrow[u, "\zresal"].\end{tikzcd}
}
\end{equation}
Note that 
$\enx_\czcc(\ma)$ is an algebra with $f\star g:=f\circ g$. Moreover,  $\hm_\czcc(\ml, \ma)$ is a left $\enx_\czcc(\ma)$-module via $f.\mu=f\circ \mu$ for any $f\in \enx_\czcc(\ma)$ and $\mu \in \hm_\czcc(\ml, \ma)$. Then the commutativity of  diagram \eqref{daj} also gives that $\cfl$ is a left $\cfc$-module via 
\beq\label{mstr}
\ch.\mu=\psi_{\ma,\unu}^{-1}(\psi_{\ml, \unu}(\ch)\circ \psi_{\ml, \unu}(\mu)),\;\text{for any}\;\ch\in \cfcc,\mu\in \cfl.
\eeq
On the other hand, since the restriction map ${\zresal}:\hm_\czcc(\ma, \ma)\ra\hm_\czcc(\ml, \ma)$ is an  $\enx_\czcc(\ma)$-module homomorphism it follows that $\resal$ is also a $\cfcc$-module homomorphism. This can be written as:
\beq\label{mstrtw}
\ch.\resal(\ch')=\resal(\ch\star\ch'),
\eeq
for any two $\ch,\ch'\in \cfcc$. In particular, for $\ch'=\epsu$ the unit of $\cfcc$ it follows that for any $\ch\in \cfcc$ one has
\beq\label{chdeps}
\ch.\resal(\epsu)=\resal(\ch).
\eeq

Define also the {\it central subspace} of $\cecc$ associated to $\ml$ as the vector space $\cel:=\hmcc(\unu, L)$. Denote by
\beq\label{celdescr}
\iotael: \cel=\hmcc(\unu, L)\hookrightarrow \hmcc(\unu, A)=\cecc
\eeq
the canonical inclusion given by $\iotael(z)=F(\iota_\ml)\circ z$ for any $z \in \cel$. 

For $z, z' \in \cel$, one can set $zz' :=m_L \circ (z \ot z')$. Then as in the case of $\cecc$ it is easy to see that  $\cel$ is a monoid with respect to this operation. Moreover, it is easy to see that $\iotael$ is a unitary algebra embedding. Since $\cecc$ is a semisimple commutative algebra over an algebraically closed field, see \cite{scalg}, it follows that $\cel$ is also a semisimple commutative algebra.

Let $\{\ell_s\}_{s\in \zir}$ be the (central) primitive idempotents of $\cel$.
One can write $\iotael_L(\ell_s)=\sum_{i \in \ccb_s}E_i$ for the decomposition of the primitive idempotents of $\cel$ inside $\cecc$. In this way we get a partition
$\zim=\ccb_0\sqcup \dots \sqcup \ccb_r.$ Through the rest of the paper, by abuse of notations we also denote the corresponding partition of the irreducible characters,  $\irr(\cc)=\ccb_0\sqcup \dots \sqcup \ccb_r$ with the same symbols $\ccb_s$, with $0\leq s\leq r$.
Without loss of generality we may  suppose  that the irreducible character $\ch_0=\epsu$ of the unit object of $\cc$ satisfies $\ch_0 \in \ccb_0$ i.e.
\beq\label{elm0}
\iotael(\ell_0)=\blam +\sum_{\{i \in \ccb_0,\;i\neq 0\}}E_i
\eeq
where $\blam$ is the idempotent integral of $\cc$. 
\br \label{pilc}
Note that the map $\pi:\ma\ra \ml$ induces also two maps at the level of characters and centers. One is a canonical embedding of vector spaces $\pilc:\cfl\hookrightarrow \cfcc,\;\mu\sent \mu\circ F(\pi).$ The other one is a surjective linear map $\pile: \cecc\ra \cel, \;z\sent \pi\circ z$.
Clearly one has $\resal\circ \pilc=\id_{\cfl}$ and $ \pile\circ \iotael=\id_{\cel}$. 
\er
There is also a nondegenerate canonical  paring given by
\beq
\langle\;,\;\rangle_\ml:\cfl \ot \cel  \ra \kk,\;\langle\alpha, z \rangle\sent \al\circ z.
\eeq 
Recall also the non-degenerate evaluation pairing $\langle,\;\rangle_{\mtr{\ma} }$ given by
\beq
\langle\;,\;\rangle_{\mtr{\ma}} :\cfcc \ot \cecc  \ra \kk,\;\langle \ch, z\rangle\sent \ch\circ z.
\eeq
\mdn
\bl 
One has the following compatibility properties between pairings:
\beq\label{comprop}
\langle \resal(\ch),\;\bar z\rangle_\ml=\langle\ch,\;{\iotael}(\bar z)\rangle_{\mtr{\ma}}.
\eeq

\beq\label{comp2}
\langle\mu, \pile(z)\rangle_\ml=\langle\pilc(\mu), z\rangle_\ma
\eeq
for any $\bar z \in \cel$, $z\in \cecc$, $\ch \in \cfcc$ and any $\mu\in \cfl$ 
\el
\bpf 
Straightforward.
\epf
\bp\label{chinbs}
For any two irreducible characters $\ch_i, \ch_j \in \irr(\cc)$ one has that  $\ch_i,\ch_j\in \ccb_s$ for some $0\leq s\leq r$ if and only if
\beqn
 \frac{\resal(\ch_i)}{d_i}=\frac{\resal(\ch_j)}{d_j}.
\eeqn
\ep
\bpf
Let ${\eta_s}\in \cfl$ be the dual basis of $\ell_s$ with respect to the non-degenerate pairing $\langle,\;\rangle_\ml$. Thus one has
$\langle \eta_t, \;\ell_s\rangle_{\mtr{L}} =\delta_{s, t}$, for any $0\leq s, t\leq r$. We will show that if $\ch_i \in \ccb_s$ then 
\beq\label{resd}
\frac{\resal(\ch_i)}{d_i}={\eta_s}.
\eeq 
This would finish the proof of the proposition.

Suppose now that $\ch_i \in \ccb_s$. Then by  Equation \eqref{comprop} one has 
$
\langle\resal(\ch_i), \ell_s\rangle_\ml=\langle\ch_i, \;{\iotael}(\ell_s)\rangle_{\mtr{\ma}} =\langle\ch_i, \;\sum_{j \in \ccb_s}E_j\rangle_{\mtr{\ma}} =d_i.$ On the other hand,  $\langle\resal(\ch_i), \ell_t\rangle_\ml=\langle\ch_i, \;{\iotael}(\ell_t)\rangle_{\mtr{\ma}} =0, \text{if}\; t \neq s.$ Therefore $\frac{1}{d_i}\resal(\ch_i)={\eta_s}$ for any irreducible character $\ch_i \in \ccb_s$. 
\epf

Recall that for the character $\ch_0$ of the unit object $\unu$ of $\cc$ we assumed that $\ch_0\in \ccb_0$ . We denote by $\epsul:=\resal(\ch_0)$ its restriction to $L$.  Note that the above proposition implies $\epsul=\eta_0$. Moreover the proof shows that $\ch_i \in\ccb_0$ if and only if 
\beq\label{cbz}
\resal(\ch_i)=d_i\epsul.
\eeq
\br
We have chosen to add the underline bar  in the notation $\epsul:L\ra \unu$ of the restriction $\resal(\epsu)$ in order to distinguished it from the counit ${\bf \eps}_L:Z(L)\ra L$ of the comonad  $Z$ evaluated at $L$.
\er
\subsection{On the right $L$-module maps in $\cfl$}
Since $\epsu \twd A\ra \unu$ is an algebra morphism in $\cc$ it follows that $\unu$ can be considered a right $A$-module in $\cc$ via the morphism $\epsu$. As explained   in \cite[Subsection 5.2]{scalg} the integral $\blam:\unu \ra A$ is the unique (up to scalar) morphism of right $A$-modules in $\cc$.

The next Lemma is an analogue of  \cite[Lemma 6.1]{scalg} and its proof follows  the same lines.
\bl \label{ez}
Let $\cc$ be a fusion category and $\ml$ be a unitary subalgebra of $\ma$ in $\czcc$.
\bne
\item
An arrow $a:\unu\ra L$ in $\cc$ is a morphism of right $L$-modules in $\cc$ if and only if for any $ z\in \cel$ one has\twd
\beq\label{modulemult}
za=\langle\epsul, z\rangle_\ml a.
\eeq
 
\item
$\ell_0$ is a morphism of right $L$-modules in $\cc$.
\ene
\el
\bpf
The proof of the first item is straightforward, see also \cite[Lemma 6.1]{scalg}.

Since $\cel$ is a product of copies of the field $\kk$ an element $\phi \in \cel$ satisfying Equation \eqref{modulemult} is unique up to a scalar. Moreover, note that the idempotent $\ell_0$ satisfies this equation. Indeed, since $\cel$ is a semisimple algebra there is $\mu\in \widehat{\cel}$ such that
$z\ell_0=\mu(z)\ell_0$ for all $z\in \cel$. Since $\langle\epsul, \ell_0\rangle_\ml=\langle \epsu, \iota^e(\ell_0)\rangle_\ma=1$ it follows that
\begin{eqnarray*}
\mu(z)&=&\langle\epsul, \mu(z)\ell_0\rangle_\ml
=\langle\epsul, z\ell_0\rangle_\ma
=\langle\epsu, \iotael(z\ell_0)\rangle_\ma
\\&=&
\langle\epsu, \iotael(z)\rangle_\ma\langle\epsu, \iotael(\ell_0)\rangle_\ma=\langle\epsu, \iotael(z)\rangle_\ma
\\&=&
\langle\epsul, z\rangle_\ml.
\end{eqnarray*}
\epf

\subsection*{On the full subcategory $\csl$}
We denote by $\csl$ the full abelian subcategory of $\cc$ generated by the simple objects $V_i$ of $\cc$ whose irreducible characters $\ch_i$ satisfy $\ch_i\in \ccb_0$.
\bt
Let $\cc$ be a pseudo-unitary fusion category and $\ml$ be a unitary subalgebra of the adjoint algebra $\ma$ of $\cc$. With the above notations one has that
${\cs_\ml}$ is a fusion subcategory of $\cc$.
\et
\bpf
Suppose that $M, M'\in \irr({\cs_\ml})$ are two simple objects with characters $\ch_M=\ch,\;\ch_{M'}=\ch'$.
Then $\resal(\ch)=d(\ch)\epsul$ and $\resal(\ch')=d(\ch)\epsul$. It follows that
\beqn
\resal(\ch\star \ch')=\ch.\resal(\ch')= d(\ch')\ch.\epsul\numeq{\ref{chdeps}} d(\ch)\resal(\ch')=d(\ch)d(\ch')\epsul
\eeqn
On the other hand suppose that
$\ch\star \ch'=\sum_{\ch_u\in \irr(\cc)}N^u_{M, M'}\ch_u$ with some $N^u_{M,M'}\geq 0$.
It follows  that
\begin{eqnarray*}
\resal(\ch\star\ch')&=&\sum_{\ch_u\in \irr(\cc)}N^u_{M, M'}\resal(\ch_u)=\sumstor \sum_{u \in \mtc B_s}N^u_{M, M'}d_u\resal(\frac{\ch_u}{d_u})\\ &\numeq{\ref{resd}}& \sumstor \big(\sum_{u \in \mtc B_s}N^u_{M, M'}d_u\big){\eta_s}
\end{eqnarray*}
where $\{\eta_s\}$ is the dual basis of $\ell_s$ with respect to the pairing $\langle,\;\rangle_\ml$. Since $\eta_s$ are linearly independent  it follows that for any $s\neq 0$ one has
$$\sum_{u \in \mtc B_s}N^u_{M, M'}d_u=0$$
Since $\cc$ is pseudo-unitary one has $d_u>0$ and therefore in this case $N^u_{M,M'}=0$ for any $u\in \mtc B_s$. Thus all the simple constituents of $M\ot M'$ are in full abelian subcategory ${\cs_\ml}$.
\epf
\subsection*{On the $\fp(\csl)$.} In this subsection we prove the following:

\bt\label{fpdsl}
Let $\cc$ be a pseudo-unitary fusion category and $\ml$ be a unitary subalgebra of the adjoint algebra $\ma$ of $\cc$. With the above notations one has that
$$\fp({\cs_\ml})=\frac{\dimcc}{\dim(\ml)}.$$ 
\et
First we need to fix several notations. Let $\cc$ be a fusion category and $F:\czcc\ra\cc$ be  the forgetful functor. As above, let $R$ be a right adjoint of $F$. We may suppose that $\ma=\bigoplus_{j=0}^r\ccpj$ is the decomposition of $\ma$ in homogenous components in $\czcc$. Define $\mtc J:=\{0,\dots, r\}$ the set of indices of homogenous components of $\ma$.

We may write each homogenous component as
$\ccpj=\bigoplus_{s=1}^{m_j}\ccpsj$
where $\ccpsj$ are the simple sub-objects of $\ma$ entering in the homogenous component $\ccpj$. Therefore as an object of $\czcc$ one has a decomposition in simple objects
\beq\label{dec}
\ma=\bigoplus_{j\in \mtc J}\bigoplus_{s\in \mtmj}\ccpjs
\eeq
where $\ccpjs$ and $\ccpjt$ are isomorphic simple $\czcc$-submodules of $\ma$ and $\mtc M_j:=\{1,\dots, m_j\}$. One has $\enx_\czcc(\ccpj)\simeq M_{m_j}(\kk)$ as algebras and therefore
$\enx_\czcc(\ma)\simeq \bigoplus_{j \in \mtc J}M_{m_j}(\kk).$
By the natural isomorphism $\phia$ from Equation \eqref{can} one also has that
$\cfcc\simeq \bigoplus_{j \in \mtc J}M_{m_j}(\kk)$
as a semisimple algebra. 

Without loss of generality we may also suppose that $\cc^{(0)}=\unu_\czcc$ is the simple trivial object of $\ma$. Since $\hm_\czcc(\ma ,\unu)\simeq \hmcc(\unu, \unu)=\kk$ it follows that $m_0=1$.

Define $\widetilde{\fjst}\in \enx_\czcc(\ma)$ as the unique endomorphisms of $\ma$ that send (as identity) $\ccjpt$ into $\ccjps$ and that are zero on the other conjugacy classes $\cc^{(j)}_{t'}$ with $t'\neq t$.  Then these endomorphisms form a linear basis $\enx_\czcc(\ma)$ such as 
$$
\wtd{F^j_{st}} \circ \wtd{F^{j'}_{s', t'}}=\delta_{j,j'}\delta_{s',t}\wtd{F^j_{st'}}
$$

Denote $\fjst:=\psi_{\ma,\unu}^{-1}(\widetilde{\fjst})\in \cfcc$. It follows that $F^j_{st}$ are the standard matrix entries in the matrix-block $M_{m_j}(\kk)$. Thus 
\beq\label{mcfa}
F^j_{st}F^{j'}_{s't'}=\delta_{j,j'}\delta_{s', t}F^j_{s t'},
\eeq
and the corresponding primitive central idempotent corresponding to the unit of the block is $F^j=\sum_{s=0}^{m_j}F^j_{ss}$. Note that the decompositions of $\ma$ from Equation \eqref{dec} are in bijection with the matrix bases $\{F^j_{st}\}$ of $\cfcc$.

It is well known that the block $\kk=M_{m_0}(\kk)$ from the above decomposition  corresponds to the central primitive idempotent  $F^0=\lam$, the cointegral of $\cc$.



By \cite[Lemma 6.2]{scalg} the Grothendieck ring $\mathrm{Gr}_\kk(\cc)$ is a symmetric Frobenius algebra with non-degenerate trace $\tau:\mathrm{Gr}_\kk(\cc)\ra \comp$ given by 
$[X]\mapsto \dimk\hm_\cc(\unu, X)$. For a pivotal fusion category one has $\tau([X])=\langle\ch(X), \blam\rangle$ for any object $X\in \cocc$, where $\blam$ is an idempotent cointegral associated to $\cc$.

Then the corresponding associative non-degenerate bilinear form on $\cfa$ is given by $\beta_\tau(\ch, \mu):=\langle\ch\mu, \blam\rangle$. A pair of dual bases for $\beta_\tau$ is given by $\{\ch_i, \ch_{i^*}\}$.

Suppose that on the matrix-block decomposition of $\cfcc$ 
one has
\beq\label{tau}
\beta_\tau(\ch,\;\mu)=\sum_j\frac{1}{n_j}\tr_j(\ch\mu).
\eeq
where $\tr_j$ is the usual trace on the mtrix algebra $M_{m_j}(\kk)$. Since  $\{\fjst,\; n_j\fjts\}$ is another pair of dual bases for $\beta_\tau(-,\;-)$ one can write:
\beq\label{teq}
\sumjtomst n_j\fjst\ot \fjts=\sumitom \ch_i\ot \ch_{i^{*}}.
\eeq

\bl 
For a pivotal fusion category, with the above notations one has
\beq\label{njclass}
n_j=\frac{\dimcc}{\dim(\cc^{(j)}_s)}.
\eeq
\el
\bpf
By \cite[Proposition 5.17]{scalg} one has that for any $f \in \cfcc$  
$$\tr(\psi_{\ma,\unu}^{-1}(f))=\langle f,\blam\rangle\langle\lam, u\rangle$$
where $(\blam, \lam)$ is any pair of a integral and cointegral of $\cc$ such that $\langle\lam,\blam\rangle=1$. One can chose such a pair with $\langle\lam,u\rangle=\dimcc$ and $\blam$ an idempotent integral, i.e $\epsu(\blam)=1$.

In particular, for $f=F^j$ one has that $\psi_{\ma,\unu}^{-1}(F^j)$ is the projection on the homogenous component $\cc^{(j)}$. Therefore the above trace formula can be written as
$$
m_j\dim(\cc^{(j)}_s)=\tr(\psi_{\ma,\unu}^{-1}(F^j))=\langle F^j, \blam\rangle\dimcc.
$$
On the other hand, from Equation \eqref{tau} one has
\beq\label{nj}
\tau(F^j_{ss})= \langle F^j_{ss}, \blam\rangle =\frac{1}{n_j}\tr_j(\fjss)=\frac{1}{n_j} 
\eeq
Thus 
$
\langle F^j, \blam\rangle=\tau(F^j)=\sum_{s=0}^{m_j}\tau(F^j_{ss})=\frac{m_j}{n_j}
$
and the result follows.
\epf
\bn{defn}
Denote by $\csu^j_{st}:=\mtf^{-1}(\fjst)\in \cecc$ and call this element the conjugacy class sum corresponding to $F^j_{st}$.  Here we use the cointegral $\lam$ with $\langle\lam, u\rangle=1$.

\end{defn}

\bl 
 Let $\cc$ be a spherical  fusion category. With the above notations one has that:
 \beq\label{otb}
\langle \fjst,\; \csu^i_{uv}\rangle_\ma =\delta_{i,j}\delta_{v,s}\delta_{u,t}{\dim(\cc^{(j)}_s)}.
\eeq
 \el
\bpf
In particular, for $\ch=\fjst$ and $\mu=\fjts$ in Equation \eqref{comp1} one has that
 \begin{eqnarray*}
\langle F^j_{st},\; \csu^i_{uv}\rangle_\ma & =& \langle F^j_{st},\;\mtf^{-1}(F^i_{uv})\rangle_\ma  \numeq{\ref{comp1}}\dimcc \tau(F^j_{st}F^i_{uv})
\\&\numeq{\ref{mcfa}}&
\delta_{i,j}\delta_{u,t}\dimcc\tau(F^i_{sv})\numeq{\ref{tau}}\delta_{i,j}\delta_{v,s}\delta_{u,t}\frac{\dimcc}{n_j}.
 \end{eqnarray*}
Then use the formula \eqref{nj}.
\epf

\bl
Let $\cc$ be a spherical non-degenerate fusion category. With the above notations one has
\begin{equation}\label{epsucj}
\langle \epsu, \csu^j_{st}\rangle_\ma =\delta_{s,t }\dim(\cc^{(j)}_s)
\end{equation}
\el
\bpf
Note that $\epsu=\ch_0=\sum_j F^j$ is the unit of $\cfcc$. Equation \eqref{comp1}  for $\ch=\epsu$ and $\mu=\fjts$ gives:
\begin{equation*}
\langle \epsu, \csu^j_{st}\rangle_\ma \numeq{\ref{comp1}}\dimcc\tau(F^j_{st})\numeq{\ref{nj}}\delta_{s,t }\frac{\dimcc}{n_j}\numeq{\ref{njclass}}\delta_{s,t }\dim(\cc^{(j)}_s).
\end{equation*}
\epf
\bl\label{flemma}
With the above notations one has
\bne
\item
$F^j_{st}\in \hmcc(\cc^{(j)}_t,\unu)$
\item
$\mtr C^j_{st}\in \hmcc(\unu, \cc^{(j)}_s)$.
\item
 On the homogenous components the Fourier transform sends
$\hmcc(\unu, F(\cc^{(j)})$ into $\hmcc(F(\cc^{(j)}),\unu)$.
\ene
\el
\bpf
By the adjunction isomorphism $\psi^{-1}_{\ma,\unu}$ from Equation \eqref{caninv} one has $F^j_{st}=\epsu\circ \widetilde{\fjst}$ which gives the proof for the first item.

On the other hand Equation \eqref{comp1} shows that also
$\{\fjst, \frac{n_j}{\dimcc}\csu^j_{ts}\}$ are also dual bases for $\langle,\;\rangle_\ma$. Since for $s\neq t$ the pairing $\langle,\;\rangle_\ma$ is zero on $\hmcc(\ccjpt, \unu)\times \hmcc(\unu,F(\ccjps))$ it follows that $\mtr C^j_{ts}\in \hmcc(\cc^{(j)}_t,\unu)$.
\epf
Recall the embedding $q^*\twd\cfcd\ra\cfcc$ from Subsection \ref{cfinc}. By abuse of notations we denote by the same symbol $\lamcd$, the image of $\lamcd\in \cfcd$ under $q^*$.

\bp\label{epsp}
Let $\cc$ be a non-degenerate spherical fusion category and $\cd\subseteq \cc$ a fusion subcategory of $\cc$.  Suppose that
$$
\lamcd=\sum_{(j,s, t)\in \mtc L_\cd}\beta^j_{st}F^j_{st}
$$ 
where all $\beta^j_{st}\in \kk$ are non-zero scalars and $\mtc L_\cd\subseteq \sqcup_{j \in \mtcj}\{j\}\times \mtc M_j\times \mtc M_j$ is a subset of indices. Then
\beq\label{eps}
\frac{\dimcc}{\dimcd}=\sum_{(j,s, s) \in \mtjcd}\beta^j_{ss}\dim(\cc^{(j)}_{s}).
\eeq
\ep

\bpf
We show that both terms of Equation \eqref{eps} equal $\langle\epsu, \ell_\cd\rangle_\ma$ where $\ell_{\cd}:=\mtf^{-1}(\lamcd)$. Here the cointegral $\lam$ is chosen such that $\langle\lam,u\rangle_\ma=1$, where $u=F(\eta_\ma)$ is the unit of $\cecc$.
Applying $\mtf^{-1}$ to the above formula for $\lamcd$ it follows that the element $\ell_{\cd}$ has the formula
\beq\label{fellcd}
\ell_{\cd}=\sum_{(j,s, t)\in \mtc L_\cd}\beta^j_{st}C^j_{st}.
\eeq
  Equation \eqref{fellcd} gives that
$$
\langle\epsu, \ell_\cd\rangle_\ma=\sum_{(j,s, t)\in \mtjcd}\beta^j_{st}\langle\epsu,\; \csujst\rangle_\ma\numeq{\ref{epsucj}}\sum_{(j,s, s)\in \mtjcd}\beta^j_{ss}\dim(\cc^{(j)}_s).
$$
On the other hand from Equation \eqref{intregform} one has
\begin{eqnarray*}
\ell_\cd&=&\mtfi(\lam_\cd)=\frac{1}{\dimcd}\big(
\sum_{\ch_i\in \irr(\cd)}d_{i^*}\mtfi(\ch_i)
\big)\\ &\numeq{\ref{610}}& \frac{1}{\dimcd}\big(
\sum_{\ch_i\in \irr(\cd)}d_{i^*}\frac{\dimcc}{d_i}E_{i^*}
\big).
\end{eqnarray*}
Since $\cc$ is spherical it follows that
\beq\label{idmptsums}
\ell_\cd=\frac{\dimcc}{\dimcd}\big(\sum_{\ch_i\in \irr(\cd)}E_{i}
\big).
\eeq
By applying $\epsu$ to Equation \eqref{idmptsums} it follows that $\langle\epsu, \ell_\cd\rangle_\ma= \frac{\dimcc}{\dimcd}$. Indeed, note that $\lag \eps_\unu, \;E_j \rag=\delta_{j,\;0}$ since we assumed that $\ch_0=\epsu$.
\epf
\br
Note that the above result holds for any decomposition of $\ma$ from  Equation \eqref{dec}, thus basically for any matrix basis $F^j_{st}$ of $\cfcc$.
\er
\br
In \cite{ccc-march} it was studied the central  element $\elcd:=\mtfi(\lam_\cd)\in \cecc$ in the case of a fusion category with commutative Grothendieck ring.
\er
To simplify the notations, through the rest of the paper, we denote $C_{st}\twd=C_{st}^j$ for the class sum $C_{st}^j$ whenever it is specified that $s, t\in \mtc M_j$, case in which the index $j$ is implicitly understood. The same notation $F_{st}\twd=F^j_{st}$ will be used in the same circumstances.

Note that the {\it central subspace} of $\cecc:=\hmcc(\unu, A)$ associated to $\ml$ can be described as $\cecc:=\bigoplus_{j \in \mtcj}\bigoplus_{s\in \mtc M_j}\hmcc(\unu,F(\ccpjs)).$ Therefore, a linear basis for $\cecc$ is given by $\csu^j_{st}$ with $j \in \mtc J$ and $s,t\in \mtc M_j$.
\subsection{On unitary subalgebras of $\ma$}\label{subalg}
Let  $\ml$ be a  subalgebra of $\ma$ in $\czcc$. Since $\ma$ is a semisimple object of $\czcc$ we may choose the simple direct summands $\cc^{(j)}_s$ from Equation \eqref{dec} such that 
\beq\label{decl}
\ml=\bigoplus_{j\in \mtcjl}\bigoplus_{s\in \lj}\ccpsj
\eeq
for some subset $\mtcjl\subseteq \mtcj$ and a subset $\lj\subseteq \mtmj:=\{1,\dots, m_j\}$. Denote by $\widetilde J_\ml:=\sqcup_{j\in \mtcjl}\{j\}\times \mtc L_j$.
\\

By the naturality of the adjunction isomorphisms $\psi_{\ccjps,\unu}$ it follows that via the inclusion $\iotael:\cel\hookrightarrow \cecc$ one has 
{
\beq\label{celdec}
\cel:=\hmcc(\unu, L)=\bigoplus_{j \in \mtcjl}\bigoplus_{s\in \lj}\hmcc(\unu,F(\ccpjs)).
\eeq
}

Therefore we may denote by the same symbol $\csu^j_{st}\in \hmcc(\unu, F(\ccjps))$ the elements of the linear basis of $\cec$ that belong to $\cel\subseteq \cecc$ if $s\in \lj$. In other words, $\iotael(\csu^j_{st})=\csu^j_{st}$ if $j\in \mtc J_\ml$ and $s \in \lj$.

Thus, a linear basis for $\cel$ is given by $\csu^j_{st}:=\mtf^{-1}(F^j_{st})$ with $j \in \mtcjl$, $s\in \lj$ and $t \in \mtc M_j$.

Note also that $\cfcc=\bigoplus_{j \in\mtcj, s\in \mtc M_j}\hmcc(F(\ccpjs), \unu)$.
Recall that there is also a projection $\pi:
\ma\ra \ml$ in $\czcc$ since this is a semisimple category. Then, as above,  the induced  linear morphism $\pilc:=F(\pi)^*:\hmcc(F(\ml),\unu)\ra\hmcc(F(\ma), \unu)$ gives an embedding $\cfl\hookrightarrow \cfcc$. Via this embedding one has  $\cfl=\bigoplus_{(j,s) \in \widetilde\mtcjl}\hmcc(F(\ccpjs), \unu)$ and we may denote by the same symbol $F^j_{st}$ the same elements corresponding to $\cfl$ when ${(j,t) \in \widetilde\mtcjl}$.

The following three lemmas are easy to verify.
\bl\label{needproofcl}
A basis for $\cfl$ is given by $\fjst$ with $(j,t) \in \widetilde\mtcjl$ and $s\in \mtc M_j$.
\el

\bl\label{resal}
With the above identifications one has
\beq\label{resalnz}
\resal(\fjst)=\fjst, \;\text{if}\;(j, t)\in \widetilde\mtcjl,\;\text{for any,}\;0\leq s\leq m_j
\eeq
 and 
\beq\label{resalz}
\resal(\fjst)=0,\;\text{if}\;(j, t)\notin \widetilde\mtcjl,\;\;\text{for any,}\;0\leq s\leq m_j.
\eeq
\el

\bl\label{piec}
For any $j \in \mtc J_L$ one has
$\pile(C^j_{s t})= C^j_{s t}$ if $j \in \mtc J_\ml$, $s\in \lj$ and, $\pile(C^j_{s t})=0$ otherwise.
\el
For any irreducible character $\ch_i$ write 
\beq\label{chi}
\ch_i=\sumjtomst \alijst F^j_{st}
\eeq
for some scalars $\alijst\in \kk$. Applying $\mtf^{-1}$ to this Equation and using  Equation \eqref{610} one has that
\beq\label{einc}
\frac{\dimcc}{{d_i} }E_{i^*}=\sumjtomst \al(i)^j_{st}\csu^j_{st}
\eeq
which gives that
\beq\label{einc2}
E_{i^*}=\frac{1}{\dimcc}\sumjtomst d_i\al(i)^j_{st}\csu^j_{st}.
\eeq
Denote by 
\beq\label{betajst}\beta^j_{st}:=\sum_{i \in \ccb_0}d_i \al(i)^j_{st}.
\eeq
\bl \label{idrange}
Let $\cc$ be a fusion category and $\ml$ be a unitary subalgebra of $\ma$. With the above notations  one has
\beq\label{idt}
\beta^j_{st}=0, \; \text{for any}\; (j,s )\notin \widetilde{\mtcjl}, \;\text{and any}\;0\leq t\leq m_j.
\eeq 
\el
\bpf Recall the linear basis $\ell_k$ of $\cel$, with $0\leq k\leq r$ .
Applying Equation \eqref{einc} one obtains that for any $0\leq k\leq r$ one has
$$
\iotael(\ell_k)=\sum_{\ch_i\in \irr(\ccb_k)}E_i\numeq{\ref{einc}}\frac{1}{{\dimcc}}\sumjtomst\bigg(\sum_{\ch_i\in \irr(\ccb_k)}{d_{i^*}\al(i^*)^j_{st}}\bigg)\csu^j_{st}.
$$
Since $\{\csu^j_{st}\}_{\jinmtcjl, s\in \lj, t\in \mtc M_j}$ is a basis of $\cel$ it follows that
\beq\label{ido}
\sum_{\ch_i\in \irr(\ccb_k)}{d_{i^*}\al(i^*)^j_{st}}=0
\eeq
if $j \notin \mtcjl$ or if $j \in \mtcjl$ but $s\notin \lj$. Note that  for $k=0$ this  can also be written as
$\beta^j_{st}=0, \; \text{for any}\; (j,s )\notin \widetilde{\mtcjl}, \;\text{and any}\;0\leq t\leq m_j
$ since $\csl=\ccb_0$ is a fusion subcategory.
\epf
\bl
With the above notations one has 
\beq
\epsul=\sum_{(j,t)\in \widetilde\mtcjl}F^j_{tt}.
\eeq
For any $\ch_i\in \ccb_0$ one has
\beq\label{lastid}
\al(i)^j_{st}=\delta_{s,t}d_i, \text{for any}\; (j,t)\in \widetilde \mtcjl,  \;\text{and any},\;1\leq s\leq m_j.
\eeq
\el
\bpf
Since $\epsu=\sumjtom F^j$ is the unit of $\cfcc$ it follows by Lemma \ref{resal} that $\epsul:=\resal(\epsu)=\sum_{(j,t)\in \widetilde\mtcjl}F^j_{tt}$. If $\ch_i\in \ccb_0$ then by Equation \eqref{cbz} one has $\resal(\ch_i)=d_i\epsul=d_i\big(\sum_{(j,t)\in \widetilde\mtcjl}F^j_{tt}\big)$. On the other hand

\begin{eqnarray*}
\resal(\ch_i)&=&\sum_{j \in \mtc J}\sum_{s\in \mtc M_j} \sum_{t\in \mtc M_j}\alijst \resal(F^j_{st})
\\&\numeq{\ref{resal}}&\sum_{j\in  \mtcjl}\sum_{ t\in \lj}\sum_{s\in \mtc M_j}\alijst F^j_{st}.
\end{eqnarray*}
Comparing the two expressions for $\resal(\ch_i)$ it follows that for $(j,t)\in \widetilde \mtcjl$ one has $\al(i)^j_{st}=d_i\delta_{s,t}$.
\epf

\bp\label{intccl}
Suppose that $\cc$ is a spherical fusion category and $\ml$ be a unitary subalgebra of $\ma$. With the above notations, for the idempotent integral of $\csl$ one  has that
\beq\label{lamcbz}
\lam_{\cs_\ml}=\sum_{j\in \mtcjl}\sum_{s \in \mtc L_j}F^j_{ss}+\sum_{j\in \mtcjl}\sum_{s \in \mtc L_j, t\notin \mtc L_j}\beta^j_{st}F^j_{st}
\eeq
for some scalars $\beta^j_{st}\in \kk$.
\ep
\bpf
Recall that $\dim({\cs_\ml})=\sum_{i\in {\irr(\csl)}}d_{i^*}d_i$. Then since $\cc$ is spherical one has that
\begin{eqnarray*}
\lam_{\cs_\ml}&=&\frac{1}{\dimcsl}\big(\sum_{\ch_i\in {\irr(\csl)}}d_{i^*}\ch_i\big)
\\&\numeq{\ref{chi}}&\frac{1}{\dimcsl}\bigg(\sum_{\ch_i\in {\irr(\csl)}}d_{i^*}\big(\sumjtom\sum_{s,t\in \mtc M_j} \al(i)^j_{st} F^j_{st}\big)\bigg)
\\ &= &\frac{1}{\dimcsl}\sumjtom\sum_{s,t\in \mtc M_j} \big(\sum_{\ch_i\in {\irr(\csl)}}d_{i}^*\al(i)^j_{st}\big)F^j_{st}
\\ &=& \frac{1}{\dimcsl}\sumjtom\sum_{s,t\in \mtc M_j} \beta^j_{st}F^j_{st}
\\ &\numeq{\ref{idt}} & \frac{1}{\dimcsl}\sum_{j\in \mtcjl}\sum_{s\in \lj}\sum_{t\in \mtc M_j} \beta^j_{st}F^j_{st}
\\&=& \frac{1}{\dimcsl}\sum_{j\in \mtcjl}\sum_{s\in \lj}\sum_{t\in \lj} \beta^j_{st}F^j_{st}
\\&+&
\frac{1}{\dimcsl}\sum_{j\in \mtcjl}\sum_{s\in \lj}\sum_{t\in \mtc M_j\setminus\lj} \beta^j_{st}F^j_{st}
\\&\numeq{\ref{lastid}}&
\frac{1}{\dimcsl}\bigg( \sum_{j\in \mtcjl}\sum_{s\in \mtc L_j} \dimcsl F^j_{ss}\bigg)\\&+&\frac{1}{\dimcsl}\sum_{j\in \mtcjl}\sum_{t\in \lj}\sum_{s\in \mtc M_j\setminus\lj}\beta^j_{st}F^j_{st}
\\ &=& \sum_{j\in \mtcjl}\sum_{t\in \mtc L_j}F^j_{ss}+\frac{1}{\dimcsl}\sum_{j\in \mtcjl}\sum_{s \in \mtc L_j, t\notin \mtc L_j}\beta^j_{st}F^j_{st}
\end{eqnarray*}
where the scalars  $\beta^j_{st}:=\big(\sum_{\ch_i\in {\irr(\csl)}}d_{i}^*\al(i)^j_{st}\big)$ are defined as in Equation \eqref{betajst}.
\epf
\subsection*{Proof of Theorem \ref{fpdsl}}

\bpf 
Note that Proposition \ref{intccl} shows that the set of indices  $\mtc L_{\csl}$ of $\csl$ coincides to $\sqcup_{j\in\mtcjl}\{j\}\times \mtc M_j\times \lj$.
Then by Equation \eqref{eps} one has that
$$\frac{\dimcc}{\dim(\csl)}=\sum_{(j,s)\in \widetilde\mtcjl}\dim(\cc^{(j)}_s)=\dim(\ml).$$
which finishes the proof of Theorem \ref{fpdsl} since  $\cc$  is a pseudo-unitary fusion category.
\epf

\bc\label{seq}
Let $\cc$ be a pseudo-unitary fusion category and $\mm$, $\ml$ be two unitary subalgebras of $\ma$.  Then $\csm=\csl$ if and only if $\mm=\ml$
\ec
\bpf
 Note that since $\ma$ is a semisimple object one can choose the   the simple direct summands $\ccjps$ 
such that Equation \eqref{decl}
is satisfied for both $\mm$ and $\ml$, i.e.
\beq\label{decl}
\ml=\bigoplus_{j\in \mtcjl}\bigoplus_{s\in \mtc L^1_j}\ccpsj,\; \mm=\bigoplus_{j\in \mtcj_\mm}\bigoplus_{s\in \mtc L^2_j}\ccpsj. \eeq
for some subsets $\mtc J_\ml, \mtc J_\mm\subseteq \mtcj$ and $\mtc L^1_j, \mtc L^2_j\subseteq \mtc M_j$. One has $\csl= \csm$ if and only if $\lam_\csl=\lam_{\csm}$. Then the result follows from Equation \eqref{lamcbz}.
\epf
\subsection{On the abelian full subcategories $\ccl$}\label{sccl}

Let $M$ be any object of $\cc$. We denote by $\racml$ the restriction to the subalgebra $L$ of the action ${\actcm}$ of $A$ on $M$ defined in Equation \eqref{rox}. We say that an object $M$ of $\cc$ {\it receives a trivial $L$-action} if
\beq\label{triviall}
\racml=\epsul\ot \id_M.
\eeq
\bn{defn} \label{ccl}
Let $\cc$ be any finite tensor category and $\ml$ be a unitary subalgebra of $\ma$.
We define ${\ccl}$ as the full abelian subcategory of $\cc$ whose objects $M$ receive a trivial action ${\actcm}$.
\end{defn}
Next proposition holds for any finite tensor category, not necessarily fusion.
\bp\label{cincls}
Let $\cc$ be a  pivotal finite tensor category and $\ml$ be a unitary subalgebra of $\ma$.
Then as abelian subcategories of $\cc$ one has  
$${\ccl}\subseteq \cs_\ml.$$ 
\ep
\bpf

Note that by definition of the dimension one has 
$$
d(M)\id_\unu=(\unu\xra{\cov_M} M\ot M^* \xra{\widetilde{ev_M}}\unu)
$$
where $\widetilde{ev_M}=(M\ot M^*\xra{j_M}M^{**}\ot M^*\xra{ev_{M^*}}\unu)$ is defined using the pivotal structure $j_M$. This implies 
\beq\label{rhs}
d(M)\epsul=(L\simeq L\ot \unu\xra{\epsul\ot \id_\unu} \unu \xra{\cov_M} M\ot M^*\xra{\widetilde{ev_M}}\unu).
\eeq

Suppose that $M\in\co(\cc)$ satisfies Equation \eqref{triviall}. We will show that $M\in \co(\csl)$. Since for any two morphisms in $\cc$ one has $f\ot g=(f\ot 1)\circ (1\ot g)$, it follows that
\beq\label{int1}
(\id_A\ot \coev_M)\circ F(\iotal)=(F(\iotal)\ot \id_M\ot \id_{M^*})\circ(\id_L\ot \coev_M).\eeq
Then the restriction $\resal(\ch(M))$ can be written as:
{
\begin{eqnarray*}
\hskip -0.7cm\resal(\ch(M))&=&(L\xra{{F(\iota_\ml)}} A\xra{\ch(M)}\unu) =
\\ &\numeq{\ref{chx}}& (L\xra{{F(\iota_\ml)}}A\simeq A\ot \unu\xra{1_A\ot \cov_M} A\ot M\ot M^*\ra
\\&\xra{{\actcm}\ot \id_{M^*}} &M\ot M^*\xra{\widetilde{\ev_M}}\unu)=
\end{eqnarray*}
\vskip -0.75cm
\begin{eqnarray*}
 &\numeq{\ref{int1}}& (L\simeq L\ot \unu \xra{\id_L\ot \cov_M} L\ot M\ot M^*\xra{\racml\ot \id_{M^*}} M\ot M^*\xra{\widetilde{\ev_M}}\unu).
\end{eqnarray*}
}
Replacing from Equation \eqref{triviall} $\racml$ by $\epsul \ot \idm$  one obtains that $\resal(\ch(M))$ can be written as
{
$$
(L\simeq L\ot \unu\xra{\id_L\ot \cov_M} L\ot M\ot M^*\xra{\epsul \ot \idm\ot \id_{M^*}} M\ot M^*\xra{\tilde{\ev_M}}\unu).
$$
}
Using Equation \eqref{rhs} this shows that $\resal(\ch(M))=d(M)\epsul$, i.e $M\in \co(\cs_\ml)$.
\epf
\section{Proof of  Theorem \ref{main1}}\label{coindm1}

Let $\cc$ be a finite tensor category and $F:\czcc\ra\cc$ be the forgetful functor. As above we denote by $R$ a right adjoint of $F$. The monoidal structure of the adjoint  functor $R:\cc\ra\czcc$ gives two actions of the algebra $\ma=R(\unu)$ on the object $R(M)\in \czcc$. The left action is defined by $\ro^l_{R(M)}=R_2(M,\unu): R(M) \ot R(\unu)\ra R(M\ot \unu)\simeq R(M)
$ and similarly the right action is defined by $\ro^r_{R(M)}=R_2(\unu, M): R(\unu) \ot R(M)\ra R(\unu\ot M)\simeq R(M).$ By \cite[Lemma 6.4]{sakalos} for any finite tensor category one has
\beq\label{sakres}
\ro^r_{R(M)}= \ro^l_{R(M)}\circ c_{R(M), \ma}
\eeq
where $c_{R(M), -}$ is the braiding of $R(M)\in \czcc$.  The author also has warned that in general 
$\ro^l_{R(M)}\neq\ro^r_{R(M)}\circ c_{\ma, R(M)}$.

\bl
Let $\cc$ be a finite tensor category. For any object $M$ of $\cc$ with the above notations one has
\beq\label{compeps}
\eps_M \circ \rho^l_{R(M)}=\epsu\ot\eps_M, \;\;\eps_M \circ \rho^r_{R(M)}=\epsu\ot\eps_M 
\eeq
\el
\bpf
It follows from the definition of the actions $\rho^l_{R(M)}, \rho^r_{R(M)}$ and Equation \eqref{ztw}.
\epf

By  \cite[Lemma 3.5]{DMNO} $\ma$ is an \'etale algebra in $\czcc$ and the right adjoint functor $R$ induces a tensor equivalence $\cc\xra{R} \czcc_{\ma}$. Recall also that  \cite[Theorem 4.10]{DMNO} gives an (inclusion reversing) bijection between the lattice of \'etale subalgebras of $\ma$ and that of fusion subcategories of $\cc$. This bijection can be described as follows:

Given an \'etale subalgebra $\ml$ of the adjoint algebra $\ma$ one considers the full subcategory $\beta(\ml)$ of $\cc$ as the category generated  by those right $\ma$-modules in $\czcc$ that are dyslectic as right $\ml$-modules with respect to  $\ml$. 
Its inverse associates to a fusion subcategory $\cd$ of $\cc$ the subalgebra $\ml_\cd=J_\cd(\unu)$ where $J_\cd$ is the right adjoint of the forgetful functor $F_\cd:\czcc\ra\czrcd.$ Recall that $\czrcd$ is the relative center of $\cd$ as a fusion subcategory of $\cc$, see \cite{DMNO}. By the proof of \cite[Theorem 4.10]{DMNO} one has
\beq\label{fpdl}
\fp(\betal)=\frac{\fp(\cc)}{\fp(\ml)}
\eeq
for any connected \'etale subalgebra $\ml$ of $\ma$.
Using Equation \eqref{sakres}, since $\ml\in \czcc$ is a sub-object of $\ma$  we note that for any $M\in \cocc$, the dyslecticity of $R(M)$ as a right $\ml$-module is  equivalent to 
\beq\label{fromsak}
\rho^l_{R(M)}\big|_{ _\ml}=\rho^r_{R(M)}\big|_{ _\ml}\circ c_{\ml, R(M)}
\eeq
For any left object $M$ of $\cc$ we denote by $\racml$ the restriction of the action  ${\actcm}:A\ot M\ra M$ of $A$ from Equation \eqref{rox} to $L$. Recall that by Equation \eqref{rox2} one has ${\actcm}=(\id_M\ot \epsu)\circ c_{\ma, M}$. Since $\ml$ is a subalgebra of $\ma$  in $\czcc$ this formula can be restricted to $\ml$ as:\beq\label{sh-brd2}
{\racml}=(\id_M\ot \epsul)\circ c_{\ml, M},
\eeq
where as above $\epsul$ denotes the restriction of $\epsu:A\ra \unu$ to $L$.

\bp \label{dysalczml}
Let $\ml$ be a unitary subalgebra of $\ma$ in $\czcc$. With the above notations if $M$ is an object of $\cc$ such that $R(M)$ is dyslectic with respect to $\ml$  then the following diagram is commutative:
{\center
\begin{tikzcd}[column sep=3cm, row sep=1cm]
 \ml\ot Z(M)
 \arrow[rr,"\alczml"] \arrow[dr, "\epsul\ot \eps_M"]  
& &
Z(M) \arrow[dl,"\eps_M"] 
  \\
& M &   
\end{tikzcd}
}

which can be written as:
\beq\label{dyscond}
\eps_M\circ \alczml=\epsul\ot \eps_M.
\eeq
\ep

\bpf
Assume that $R(M)$ is a dyslectic right  $\ml$-module in $\czcc$. From Equation \eqref{fromsak} it follows that 
\beq
\eps_M \circ \rho^l_{R(M)}\big|_{ _\ml}=\eps_M \circ \rho^r_{\zom}\circ \rho^r_{R(M)}\big|_{ _\ml}\circ c_{\ml, R(M)}.
\eeq
By Equation \eqref{compeps} the above equality can be written as
\beq\label{fromsak2}
\epsul\ot\eps_M= (\epsul\ot\eps_M)\circ c_{\ml, R(M)}.
\eeq

On the other hand by Equation \eqref{sh-brd2} one has
\beqn
(\eps_M\ot\epsul)\circ c_{\ml, R(M)}=(\eps_M\ot \id_\unu)\circ [(\id_{R(M)}\ot \epsul)\circ c_{\ml, Z(M)}]=(\eps_M\ot \id_\unu)\circ {\alczml}.
\eeqn
Thus Equation \eqref{fromsak2} can be written as 
$\epsul\ot \epsm=\eps_M\circ \alczml$.
\epf
\subsection{Proof of the inclusion $\betal\subseteq {\ccl}$}
For any unitary subalgebra $\ml$ of $\ma$ we denote by $\beta(\ml)$ the full subcategory of $\cc$ which by the tensor equivalence $R\twd\cc\ra\czcc_\ma$ correspond to the right $\ma$-modules that are dyslectic with respect to $\ml$. Note that as mentioned above,  if $\ml$ is an \'etale subalgebra then $\beta(\ml)$ is a fusion subcategory of $\cc$ by \cite[Theorem 4.10]{DMNO}.
\bp
\label{betainccl}
Let $\cc$ be a fusion category and $\ml$ be a unitary subalgebra of $A$. With the above notations one has  $\betal\subseteq {\ccl}$.
\ep
\bpf
Suppose that $M$ is an object of $\betal$, i.e a dyslectic right $\ml$-module. By the previous Proposition the action of $\ml$ on $Z(M)$ satisfies $\eps_M\circ\alczml=\epsul\ot \eps_M$.

Since $\eps_M:Z(M)\ra M$ is a morphism in $\cc$ it follows that it  is also  an $A$-linear morphism, see \cite{scalg}, and by restriction it is also $L$-linear. 
Therefore we have
\beq\label{alinear}
\eps_M\circ \alczml=\alcml \circ (\id_\ml\ot \epsm)
\eeq
Then Equation \ref{dyscond} can be written as
\beq\label{fourstar}
(\ml\ot Z(M)\xra{\epsul\ot \epsm}M)=(\ml\ot Z(M)\xra{\id_\ml \ot \epsm} \ml\ot M\xra{{\actcm}}M).
\eeq
Then \cite[Proposition 5.1]{brn} applied to $\cc^{\op}$ implies that the counit $\eps_M$ is an epimorphism.  Since $\cc$ is semisimple it follows that $\epsm:Z(M)\ra M$ has a section $r_M:M\ra Z(M)$ in $\cc$, i.e. a morphism satisfying $\eps_M\circ r_M=\id_M$. Composing Equation \eqref{fourstar} with $(\id_\ml \ot r_M)$ one obtains that
\begin{eqnarray*}
& & (\ml\ot M\xra{\id_\ml\ot r_M}\ml\ot Z(M)\xra{\epsul\ot \epsm}M)=\\&=&(\ml\ot M\xra{\id_\ml\ot r_M}\ml\ot Z(M)\xra{\id_\ml \ot \epsm} \ml\ot M\xra{\al_M}M)
\end{eqnarray*}
which can be written as
\beqn
(\ml\ot M\xra{\epsul\ot \id_M}M)=(\ml\ot M\xra{{\actcm}}M).
\eeqn
This shows that $M\in {\ccl}$.
\epf
\subsection{Proof of Theorem \ref{main1}}
\bpf
First, we show that 
 $\beta(\mm)= \cc_\mm^{\mtr{triv}}=\cs_{\mm}$
 for any \'etale subalgebra $\mm$ of  $\ma$. Indeed, Proposition \ref{cincls} and Proposition \ref{betainccl} imply that 
$\beta(\mm)\subseteq \cc^{\mtr{triv}}_\mm\subseteq \cs_{\mm}$. Moreover by Theorem \ref{fpdsl}  the two fusion subcategories $\beta({\mm})$ and $\cs_{\mm}$  have the same Frobenius-Perron dimensions, namely $\frac{\fpcc}{\fp(\mm)}$.

On the other hand for any unitary subalgebra $\ml$ of $\ma$ we know that $\cs_\ml$ is a fusion subcategory of $\cc$ with dimension $\frac{\fpcc}{\fp(\ml)}$. The bijective correspondence $\beta$ from \cite{DMNO} gives that   $\csl=\betam$ for some \'etale subalgebra $\ml$ of $\ma$. The first part of the proof  implies that $\csl=\cs_\mm$ and Corollary \ref{seq} implies that $\mm=\ml$. 
\epf
\subsection{Proof of Theorem \ref{main2}}
\newcommand{\repca}{{\rep_\cc(A)}}
Let $\cc$ be a braided fusion category  $F:\cc\ra\cd$ a tensor functor to an arbitrary fusion category $\cd$. Recall that a {\it central structure} for $F$ is a braided  functor $\cc\xra{F'}\cz(\cd)$ such that the composition $\cc\xra{F'}\cz(\cd)\xra{\mtr{Forg}}\cd$ coincides to $F$. A functor $F$ admitting such a structure is called central. It follows by \cite[Lemma 3.5]{DMNO} that $R(\unu)$ is an \'etale algebra in $\cc$ where $R$ is a right adjoint $R$ of the central functor $F$. 
\bpf
By \cite[Lemma 3.29]{DMNO} one has that the free functor $F_A:\cc\ra\cc_A$ is a central functor and moreover, $A=R_A(\unu)$ where $R_A$ is the right adjoint of $F_A$. Since $\cc$ is non-degenerate, there is an equivalence of categories $\czcc \simeq \cc\bxt \cc^{\rev}$. By \cite[Corollary 3.30]{DMNO}  there is also an  equivalence of categories 
$\cz(\cc_A)\simeq \cc\bxt (\cc^0_A)^{\rev}.$ Moreover, under this equivalence, the free functor $F_A$ identifies as
$$
(F_A:\cc\ra \cc_A)=(\cc=\cc\bxt\vect\xra{\iota}\cc\bxt (\cc^0_A)^{\rev}\simeq \cz(\cc_A)\xra{\mtr{Forg}} \cc_A).
$$
which shows that $R_A=\iota^{\mtr{ra}}\circ R$ where $R$ is the right adjoint functor of $\mtr{Forg}:\cz(\cc_A)\ra\cc_A$ and $\iota^{\mtr{ra}}$ is the right adjoint of the inclusion $\iota:\cc\hookrightarrow \cc\bxt (\cc^0_A)^{\rev}$. Then $A=R_A(\unu)=\iota^{\mtr{ra}}(R(\unu))$ is a subalgebra of $R(\unu)$ as being the largest sub-object of $\cz(\cc_A)$ that belongs to $\cc$. Then any unitary subalgebra of $A$ in $\cc$ is also a  subalgebra of $R(\unu)$ in $\cz(\cc_A)$ and therefore \'etale by Theorem \ref{main1}. 
\epf
\br
It is interesting to investigate if the result of Theorem \ref{main2} can  be deduced by directly investigating the semisimplicity of the corresponding module category $\cc_M$. However, note that our proof gives a new perspective on the bijective correspondence from \cite{DMNO}. It shows that at least in the pseudo-unitary case, analogue to the group representations situation, fusion subcategories can be obtained via trivial actions or corresponding kernels of objects. Another natural question at this stage is if the above result holds without the pseudo-unitarity assumption on the category or if one can remove
the non-degeneracy assumption in the above statement.
\er
\subsection{A diagonal formula for $\lamcd$}\label{aptwo}
In this subsection we assume that $\cc$ is pseudo-unitary. For  a unitary subalgebra $\ml$ of $\ma$ in $\czcc$ we let $\cd=\csl$. 
Clearly, the algebra inclusion $\iota:\ml\hookrightarrow \ma$ is an inclusion in the module category $\czcc_\ml$. Since $\ml$ is an \'etale algebra it follows by \cite[Theorem 3.2]{DMNO}  that $\czcc\;_\ml$ is a semisimple category. The semisimplicity of this category implies that $\iota$ admits a retract $\pi:\ml\ra \ma$ in the same category. For the rest of this section we assume that $\pi$ is such a retract.

It is easy to check that
$\mtf(\blam)=\frac{1}{\dimcc}\epsu.$ Since $\epsu=\sum_{j\in \mtc J} \sum_{ s\in \mtc M_j}F^j_{ss}$ follows that
\beq\label{blamclss}
\dimcc \blam=\sum_{j\in \mtc J, s\in \mtc M_j}C^j_{ss}.
\eeq
Applying $\pile$ to the Equation \eqref{blamclss} it follows that
\beq\label{pieblami}
 \pile(\blam)=\frac{1}{\dimcc}\big(\sum_{j \in \mtc J_L}\sum_{s\in \lj}C^j_{s s}\big)
\eeq

\bl
With the above chosen map $\pi$ one has
\beq
\pile(\blam)=\frac{1}{\dimcd} \ell_0.
\eeq
\el
\bpf
Since by construction  $\pile(\blam)\twd=\unu \xra{\blam} A\xra{\pi} L$ is  a morphism of right $L$-modules in $\cc$ it follows by Lemma \ref{ez} that there is $\al \in \kk$ such that $\pile(\blam)=\al \ell_0.$ Equation \eqref{pieblami} implies that 
$
\dimcc \al \ell_0=\sum_{j \in \mtc J_L}\sum_\sinl C_{\bls \bls}
$
Applying $\epsul$ to this equation one has
$\dimcc \al=\diml$ which implies that $\al=\frac{\diml}{\dimcc}=\frac{1}{\dimcd}$.
\epf

\br
 Equation \eqref{pieblami} implies now that
\beq\label{ezdiag}
\ell_0={\dimcd}\pie(\blam)=\frac{\dimcd}{\dimcc}(\sum_{j \in \mtc J_L}\sum_\sinl C_{\bls \bls})
\eeq
\er

\bp
With the above notations it follows that 
\beq\label{lamcddiag}
\lam_{\cd}=\sum_{j\in \mtcjl}\sum_{s \in \mtc L_j}F^j_{ss}.
\eeq
\ep
\bpf
Since $\ellcd:=\mtfinv(\lam_\cd)$ one has by Equation \eqref{lamcbz} that
\beq\label{nfellcd}
\ellcd=\sum_{j\in \mtcjl}\sum_{s \in \mtc L_j}C^j_{ss}+\frac{1}{\dim(\csl)}\big(\sum_{j\in \mtcjl}\sum_{s \in \mtc L_j, t\notin \mtc L_j}\beta^j_{st}C^j_{st}\big).
\eeq
Using Lemma \ref{piec} it follows that
\beq\label{piewellcd}
\pie(\ellcd)=\sum_{j\in \mtcjl}\sum_{s \in \mtc L_j}C^j_{ss}+\frac{1}{\dim(\csl)}\big(\sum_{j\in \mtcjl}\sum_{s \in \mtc L_j, t\notin \mtc L_j}\beta^j_{st}C^j_{st}\big).
\eeq
On the other hand $\ell_0=\pie(\iotael(\ell_0))=\pie(\sum_{i\in \ccb_0}E_i)=\frac{\dimcd}{\dimcc}\pie({\ell_\cd}).$ Comparing Equations \eqref{ezdiag} and \eqref{piewellcd} it follows that $\beta^j_{st}=0$ for $s\neq t$.
\epf

\subsection{On the lattice of fusion subcategories} In this subsection we prove Theorem \ref{main3}. It follows from the next three results presented below. For two subalgebras of the adjoint algebra $\ma$ we denote by $\ml\mm$ the image of $\ml\ot \mm$ under the multiplication $m:\ma\ot \ma\ra \ma$. Since $\ma$ is commutative it follows that $\ml\mm =\mm\ml$ and moreover $\mm \ml$ is a subalgebra of $\ma$. Note that $\ml\mm $ is the smallest unitary subalgebra of $\ma$ containing both $\ml$ and $\mm$.

\bp
Let $\cc$ be a pseudo unitary fusion category. With the above notations one has
$$\csl\cap\csm=\cs_{\ml\mm}.$$
\ep
\bpf 
Corollary \ref{seq} implies that $\cs_{\ml\mm}\subseteq \cs_{\ml}\cap \cs_{\mm}$ since $\ml, \mm\subseteq \ma$. On the other hand, if $V\in \cs_{\ml}\cap \cs_{\mm}$, then $\ml\mm$ acts trivially on $V$, by the associativity of the action $\al_V$ of $A$ on $V$.
This shows that $\cs_{\ml\mm}\supseteq \cs_{\ml}\cap \cs_{\mm}$. 
\epf
\bl 
Let $\cc$ be a  pseudo-unitary fusion category and $\ml,\mm$ be two unitary subalgebras of $\ma$.  Then:
\beq\label{ineq}
\fp(\csl\vee \csm)\geq \fp(\cs_{\ml\cap\mm})
\eeq
\el
\bpf
Suppose ${\csl\vee \csm}=\cs_\mbp$ for some subalgebra $\mbp\subseteq \ma$.
One has $\csl\subseteq \cs_\mbp$ which implies that $\mbp\subseteq \ml$ since $\beta$ is a reversing order bijection. Similarly one deduces that $\mbp\subseteq \mm$, and therefore $P\subseteq \ml\cap \mm$.

Since $\ma$ is a semisimple object, similarly to Corollary \ref{seq} one can choose a decomposition of $\ma$ from equation \eqref{dec} such that 
$$\ml=\bigoplus_{j\in \mtc J_\ml}\bigoplus_{s\in \mtc L^1_j}\cc^{(j)}_s,\;\mm=\bigoplus_{j\in \mtc J_\mm}\bigoplus_{s\in \mtc L^2_j}\cc^{(j)}_s,\;\mbp=\bigoplus_{j\in \mtc J_\mbp}\bigoplus_{s\in  \mtc{L}^3_j}\cc^{(j)}_s.$$
for some subsets $\mtc L^i_j\subseteq\mtc M_j=\{0,1,\dots ,m_j\}$ for $i=1,2,3$.
Note also that
$$\ml\cap \mm=\bigoplus_{j\in \mtc J_\ml \cap \mtc J_\mm}\bigoplus_{s\in  \mtc{L}^1_j\cap \mtc L^2_j}\cc^{(j)}_s.$$
Since $\mbp\subseteq \ml\cap \mm$ it follows that $\mtc J_\mbp\subseteq \mtc J_\ml\cap \mtc J_\mm$ and $\mtc L_j^3\subseteq \mtc L_j^1\cap \mtc L_j^2$ for any $j \in \mtc J_\mbp$.

By Equation \eqref{lamcddiag} one has 

\beq\label{ld2}
\lam_{\csl}=\sum_{j\in \mtcjl}\sum_{s \in \mtc L^1_j}F^j_{ss},\; \lam_{\csm}=\sum_{j\in \mtc J_\mm}\sum_{s \in \mtc L^2_j}F^j_{ss},\;\lam_{\cs_{\ml\cap \mm}}=\sum_{j\in \mtc J_\ml \cap \mtc J_\mm}\sum_{s\in  \mtc{L}^1_j\cap \mtc L^2_j}F^j_{ss}.
\eeq
Also
$$
\lam_{\cs_{\mbp}}=\sum_{j\in \mtc J_\mbp }\sum_{s\in  \mtc{L}^3_j\cap \mtc L^2_j}F^j_{ss}.
$$
By Equation \eqref{nj} and Equation \eqref{njclass} one has 

$$
\tau(F^j_{ss})=\frac{1}{n_j}=\frac{\dim(\cc^{(j)}_s)}{\dimcc} >0
$$
since $\cc$ is pseudo-unitary.
It follows that 
$$
\tau(\lam_{\cs_\mbp})=\sum_{j \in \mtcj_\mbp}\sum_{s\in \mtc L^3_j }\tau(F^j_{ss})=\sum_{j \in \mtcj_\mbp}\frac{|\mtc L^3_j|}{n_j}.
$$
Similarly one deduces that 
$$\tau(\lam_{\cs_{\ml\cap \mm}})=\sum_{j \in \mtcj_\ml\cap \mtcj_\mm}\frac{|\mtc L^1_j\cap \mtc L^2_j|}{n_j}.
$$
Since  $\mtc J_\mbp\subseteq \mtc J_\ml\cap \mtc J_\mm$ and for any $j \in \mtc J_\mbp$ one has $\mtc L_j^3\subseteq \mtc L_j^1\cap \mtc L_j^2$ 
it follows that
\beq\label{ineq1}
\tau(\lam_{\cs_\mbp})\leq \tau(\lam_{\cs_{\ml\cap\mm}}).
\eeq
On the other hand, since
$
\lam_{\cs_\mbp}=\frac{1}{\dim({\cs_\mbp})}\big(\sum_{\ch_i\in \irr(\cs_\mbp)}d_i\ch_i\big)
$
we have that $\tau(\lam_{\cs_\mbp})=\frac{1}{\dim({\cs_\mbp})}$ since $\tau(\ch_i)=0$ if $\ch_i\neq \epsu$.
Similarly
$\tau(\lam_{\cs_{\ml\cap \mm}})=\frac{1}{\dim({\cs_{\ml\cap \mm}})}$. 
Then Equation \eqref{ineq1} can be written as 
$$\frac{1}{\dim({\cs_{\ml\cap\mm}})}\geq \frac{1}{\dim(\cs_\mbp)}.$$
which proves the lemma.
\epf

\bp
Let $\cc$ be a pseudo unitary fusion category. With the above notations one has
$$\csl\vee \csm=\cs_{\ml\cap \mm}$$
\ep
\bpf
Since $\csl\subseteq\cs_{\ml\cap \mm}$ and $\csm\subseteq\cs_{\ml\cap \mm}$ it follows that $\csl\vee \csm\subseteq  \cs_{\ml\cap \mm}$. The equality follows from the Lemma above.
\epf
\bp
With the above notations one has 
$$\fp(LM)\leq\frac{\fp(\ml)\fp(\mm)}{\fp(\ml\cap \mm)}.$$ Moreover one has equality if and only if $\csl\csm=\csm\csl$.
\ep
\bpf
From \cite[Lemma 3.2]{g-fact} one has
$$
\fp(\csl\csm)=\frac{\fp(\csl)\fp(\csm)}{\fp(\csl\cap \csm)}
$$
which gives that
$$
\fp(\csl\csm)=\frac{\fp(\cc)\fp(\ml\mm)}{\fp(\ml)\fp(\mm)}
$$
Since $\csl \csm\subseteq \csl\vee \csm$, one has 
$\fp(\csl \csm)\leq \fp(\csl\vee \csm)$. Note that
$$\fp(\csl\vee \csm)=\fp(\cs_{\ml\cap \mm})=\frac{\fp(\cc)}{\fp(\ml\cap\mm)}.$$ Thus
$$
\frac{\fp(\cc)}{\fp(\ml\cap\mm)}\geq \frac{\fp(\cc)\fp(\ml\mm)}{\fp(\ml)\fp(\mm)}
$$
which gives the above inequality. One has equality above if and only if $\csl\csm=\csl\vee\csm$ which is equivalent to $\csl\csm=\csm\csl$.
\epf
\bc 
If $\cc$ is a pseudo-unitary fusion category with a commutative grothendieck ring then:
One has that
$$\fp(\ml\mm)=\frac{\fp(\ml)\fp(\mm)}{\fp(\ml\cap \mm)}$$
for any tow uniray subalgebras $\ml$ and $\mm$ of $\ma$.
\ec
\bpf
It is like Maschke's theorem:
\epf

\subsection{Examples from semisimple Hopf algebras}\label{ha}
Let $H$ be a semisimple Hopf algebra and let $\cc=\rep(H)$ be the fusion category of its finite dimensional representations over $\kk$.  As in \cite[Section 3.7]{scalg} we identify the left monoidal center  of $\cc$ with the category $\czcc\simeq ^H_H\mtc{YD}$ of left-left Yetter-Drinfeld modules of $H$. The objects of the former category are pairs $(V, \delta)$ where $V$ is a left $H$-modules and $\delta:V\ra H\ot V$ defines a left $H$-comodules structure on $V$ such that
$$
(hv)\dmino\ot (hv)\dzero=h_1v\dmino S(h_3)\ot h_2v\dzero
$$
for any $h\in H$ and any $v\in V$. A morphism in the category $^H_H\mtc{YD}$ is a linear map which is simultaneously a morphism of left $H$-modules and of left $H$-comodules.  The left half-braiding of $V\in\;^H_H\mtc{YD}$ can be written as
$$c_{V, M}:V\ot M\ra M\ot V,\; v\ot m\sent v_{-1}m\ot v_0$$
for any object $M\in \rep(H)$.
The forgetful functor $$F:^H_H\mtc{YD}\ra H-\text{mod}$$
admits a right adjoint $$R:H-\text{mod}\ra ^H_H\mtc{YD}$$ which can be written as $R(M)=H\ot M$ with $$h.(l\ot m)=h_1lS(h_3)\ot h_2m,\;\delta(h\ot m)=h_1\ot (h_2\ot m).$$ Thus, with the above notations, the Hopf comonad $Z:=FR$ can be written as $Z(M)=F(R(M))=H\ot M$ with the left $H$-module structure defined above.  Then the left half-braiding of $R(M)$ becomes:
$$
c_{\nzm, X}:\nzm\ot X\ra X\ot \nzm,\; (h\ot m)\ot x\sent h_1x\ot (h_2\ot m)
$$
for any $X\in \rep(H)$.
In particular one has that $\ma:=R (\unu)$ is the pair $(H_{\ad}, \Delta)$ where $H_{\ad}$  coincides to $H$ as vector spaces, and  the  structure as left $H$-module on $H_{\ad}$ is given by  the left  adjoint action $h.a=h_1aS(h_2)$.  Moreover the left comodule structure on $H_{\ad}$ is given by the comultiplication $\bdelta: H_{\ad}\ra H\ot H_{\ad}$. The multiplication and the unit of $R(\unu)$ are given the usual multiplication and unit of $H$. The universal dinatural projections from Subsection \ref{hcm} are given by: $$\pi_{\unu;\;M}:H_{\ad}\ra  M\ot M^*,\;\; h\mapsto \sum_i  hm_i\ot m_i^*,$$
where $m_i$ is a $\kk$-linear basis on $M$ and $\{m_i^*\}$ is the corresponding dual basis in $M^*$.
In this case the natural action of $A=H_{\ad}$ from Equation \eqref{rox}  can be written as  $\al_M:A\ot M\ra M, \al_M(a\ot m)=a.m$, for each object $M\in \rep(H)$. Moreover, by considering the standard pivotal structure on $\rep(H)$, it follows that  $\cfcc=C(H)$, the character ring of $H$, and $\cecc=Z(H)$, the center of $H$. 

In this case, unitary subalgebras of $\ma\in \czcc$ are in bijection with the left normal coideal subalgebras of $H$.
Recall that $\ml\subseteq H$ is a  left normal coideal subalgebra of $H$ if it is a subalgebra of $H$ with $\bdelta(\ml)\subseteq H\ot \ml$ and $\ml$ is closed under the left adjoint action. 

Moreover, for a semisimple Hopf algebra $H$  it follows that $\beta(\ml)=\csl=\ccl=\rep(H//\ml)$ is the fusion subcategory of $\rep(H)$ consisting of those $H$-modules receiving trivial action from $\ml$. The above equalities in this case were previously shown in \cite{iop}.

The character space $\cfl:=\hm_H(L,\unu)$ defined in Section \ref{charthry}  it can be described by $\mu\in \cfl\iff \mu(h_1lS(h_2))=\eps(h)\mu(l)$, for all $h \in H$ and $l \in L$. Note that there is an inclusion of  $$\cfl:=\hm_H(L,\unu) \xhookrightarrow{\nu} \hm_L(L,\unu)=C(L),$$ into the character space (trace space)  $C(L)$ of left $L$-modules.

Denote by $\res^A_L\twd H-\mtr{mod}\ra L-\mtr{mod}$ the restriction functor of $H$-modules to the subalgebra $L$. This functor induces a map ${(\res^A_L)}_*\twd \cfcc \ra C(L)$ at the level of characters. It is easy to see that the image of this map is $\cfl\xhookrightarrow{\nu} C(L)$. 

Note also that the central space $\cel$ becomes in this case $\cel=\hm_H(\unu, L)=L\cap Z(H)$ and the inclusion $\iotael:\cel\ra \cecc$ is the usual space inclusion $L\cap Z(H)\hookrightarrow Z(H)$.
Note also that in this case $\ell_\cd$ is a multiple of the integral $\blam_L$ defined in \cite{gmj}.
\bibliographystyle{plain}
\bibliography{ccts}

\ed
\newpage
\subsection*{New directions}-\\
Look at the dimension of $\cfcc \lamcd$ and study this space.

Don't forget ccts-braided from comp B.

\blue{\Small
\br
$\overline{\al(i)^j_{st}}=\al(i^*)^j_{st}$and $\overline{\dim(V_i)}=\dim(V_{i^*})$. Do I need this remark?
\er
}
\section*{Rest of results}
\bl
On the other hand we will show that 
\beq\label{shgen}
\langle\ch_r, \csu^i_{uv}\rangle_\ma=\al(r)^i_{vu}{\dim(\cc^{(i)}_v)}
\eeq
\el
\bpf
Indeed,
\beqn
\langle\ch_r, \csu^i_{uv}\rangle_\ma\numeq{\ref{chi}}\sumjtomst \al(r)^j\st \langleF^j_{st}, \csu^i_{uv}\rangle_\ma\numeq{\ref{otb}}\sum_{s,t\in \mtc M_i}\delta_{v,s}\delta_{u,t}\al(r)^i\st{\dim(\cc^{(i)}_s)}
\eeqn
$$
=\al(r)^i_{vu}{\dim(\cc^{(j)}_v)}.
$$
\epf
\bl \label{511}
Let $\cc$ be a  non-degenerate spherical fusion category. With the above notations one has that:
\beq\label{cjc}
\csu^i_{uv}=\frac{\dimcc}{n_i}\big(\sumrtom \frac{\al(r)^i_{vu}}{d_r} E_r\big)
\eeq
\el
\bpf
Since $\{\ch_i,\;\frac{ E_i}{d_i}\}$ are dual bases for $\langle, \;\rangle_\ma$ one has 
$$
\csu^i_{uv}=\sumrtom  \langle\ch_r, \csu^i_{uv}\rangle_\ma \frac{E_r}{d_r}\numeq{\ref{shgen}}{\dim(\cc^{(i)}_v)}\big(\sumrtom\al(r)^i_{vu}\frac{E_r}{d_r}\big)
$$
and then use Equation \eqref{njclass}.
\epf
\subsection*{Characterization of $\csl$}-\\
Denote by $d_M$ the quantum dimension of $M$.
\bp
One has that $\ch_M\in {\cs_\ml}$ if and only if $\ch_M\star F^j_{ss}=d_MF^j_{ss}$ for any $s\in \lj$.
\ep
\bpf
One has that
\beq\label{eq1}
F^\jp_{\spr,\tp}F^j_{st}\numeq{\ref{mcfa}}\delta_{t,\spr}F^j_{s\tp}
\eeq
If $M\in {\cs_\ml}$ then based on the above Equation\eqref{mcfa} one has
$$
\ch_M\star F^j_{ss}=\sum_{j,\spr,\tp}\al(M)^j_{\spr,\tp}F^j_{\spr,\tp}\star F^j_{ss}\numeq{\ref{eq1}}\sum_{j,\tp}\al(M)^j_{st'}F^j_{st'}\numeq{\ref{slco}}\al(M)^j_{ss}F^j_{ss}.
$$
\epf
\blue{\Small
Do I need this proposition in the proof of $\lam_\csl$?
}
\bp 
Suppose that $\cc$ is a \blue{non-degenerate}  spherical fusion category and $\lam$ is an integral with $\lag\lam, u\rag=1$. With the above notations:
\beq\label{idmptsums}
\elcd=\frac{\dimcc}{\dimcd}(\sum_{[V_i]\in \irr(\cd)}E_i).
\eeq
\ep
\bpf
One has by definition
$
\lam_\cd=\frac{1}{\dimcd}\big(
\sum_{\ch_i\in \irr(\cd)}d_{i^*}\ch_i\big)
$
and therefore
\begin{eqnarray*}
\ell_\cd&=&\mtfi(\lam_\cd)=\frac{1}{\dimcd}\big(
\sum_{\ch_i\in \irr(\cd)}d_{i^*}\mtfi(\ch_i)
\big)\\ &\numeq{\ref{610}}& \frac{1}{\dimcd}\big(
\sum_{\ch_i\in \irr(\cd)}d_{i^*}\frac{\dimcc}{d_i}E_{i^*}
\big)
\end{eqnarray*}
The result follows since $\cc$ is spherical.
\epf

\blue{\Small
Need to say the above proposition works word-by -wors as in Prop ?? \cite{march-ccc}.
}

\blue{\Small 
Define the subset $\mtjcd$ as the set of triples such that $\beta^j_{st}\neq 0$.
}

\section*{Was it wrong?}
From the adjunction of $F$ and $R$ one has
\beq\label{dimhmccj}
\hmcc(F(\ccjps),\unu)\simeq \hm_\czcc(\ccpjs, \ma)\simeq \kk^{m_j}.
\eeq
Since $\cc$ is a fusion category \blue{or maybe because it is modular, see \cite{scalg}} it also follows that
\beq\label{dimhmccjcf}
\hmcc(\unu, F(\ccjps))\simeq \hm_\czcc( \ma, \ccpjs)\simeq \kk^{m_j}.
\eeq
\green{
\bl\label{flemma}
Need to show that Fourier transform sends
$$
\mtf:\hmcc(\unu, F(\ccjps))\ra \hmcc(F(\ccjps),\unu)
$$
\el
\bpf
Define  $\beta^l:\cfcc\ra\cecc^*$ by 
$$
\beta^l(\ch)(z)=\langle\ch, \;z\rangle_\ma
$$ for any $\ch \in \cfcc$ and $z\in \cecc$.
Since $\langle\ch, \;z\rangle_\ma=\tr_q(\ch\circ z)$ it follows that $\beta^l(\hmcc(\fcjps,\unu)\in \cecc^*$ has non-zero values only on the space $\hmcc(\unu, \fcjps)$.
\\
Moreover, Equation \eqref{comp1} shows that $\csu^i_{st}:={\mtf}^{-1}(F^i_{st})\in \hmcc(\unu, \fcjps)$.
\epf
}
\blue{\Small
Moreover, Equation \eqref{comp1} shows that $\csu^i_{ts}:={\mtf}^{-1}(F^i_{ts})\in \hmcc(\unu, \fcjps)$.
}
\bpf
The proof follows from Equation \eqref{comp1} of the two forms. By definition $F^j_{st}\in \hmcc(\fcjps, \unu)$ and $\beta^l(\fjst):\cecc\ra \kk$ sends $\hmcc( \unu, \fcjpt)$ to zero for any $t\neq s$. 
\epf
\green{
\bpf
Define  $\beta^l:\cfcc\ra\cecc^*$ by 
$$
\beta^l(\ch)(z)=\langle\ch, \;z\rangle_\ma
$$ for any $\ch \in \cfcc$ and $z\in \cecc$.
Since $\langle\ch, \;z\rangle_\ma=\tr_q(\ch\circ z)$ it follows that $\beta^l(\hmcc(\fcjps,\unu)\in \cecc^*$ has non-zero values only on the space $\hmcc(\unu, \fcjps)$.
\\
Moreover, Equation \eqref{comp1} shows that $\csu^i_{st}:={\mtf}^{-1}(F^i_{st})\in \hmcc(\unu, \fcjps)$.
\epf
}
\blue{\Small
Moreover, Equation \eqref{comp1} shows that $\csu^i_{ts}:={\mtf}^{-1}(F^i_{ts})\in \hmcc(\unu, \fcjps)$.
}
\bpf
The proof follows from Equation \eqref{comp1} of the two forms. By definition $F^j_{st}\in \hmcc(\fcjps, \unu)$ and $\beta^l(\fjst):\cecc\ra \kk$ sends $\hmcc( \unu, \fcjpt)$ to zero for any $t\neq s$. 
\epf

Maybe it has a greater impact generalizing cw for dual bases and ccc-march from monat.
\bibliographystyle{alpha}
\bibliography{24nov}
\ed